\documentclass[]{amsart}
\usepackage[T1]{fontenc}
\usepackage{graphicx}
\usepackage{amssymb}
\usepackage[all]{xy}

\makeatletter
  \theoremstyle{plain}
  \newtheorem{thm}{Theorem}[section]
  \numberwithin{equation}{section} 
  \numberwithin{figure}{section} 
  \theoremstyle{plain}
  \newtheorem{lem}[thm]{Lemma} 
  \theoremstyle{plain}
  \newtheorem{cor}[thm]{Corollary} 
  \theoremstyle{plain}
  \newtheorem{prop}[thm]{Proposition} 
  \theoremstyle{definition}
  \newtheorem{rem}{Remark}[section]
  \theoremstyle{remark}
  \newtheorem*{rem*}{Remark}

\usepackage{psfrag}
  \usepackage{bbm}

\DeclareMathOperator{\card}{card}

\renewcommand{\Re}{\mathfrak{Re}}
\renewcommand{\Im}{\mathfrak{Im}}

\newcommand{\Z}{\mathbb{Z}}
\newcommand{\R}{\mathbb {R}}

\renewcommand{\P}{\mathcal{P}}
\newcommand{\PP}{\mathfrak{P}}

\renewcommand{\H}{\mathbb {H}}

\newcommand{\N}{\mathbb {N}}

\newcommand{\I}{\mathbb {I}}
\renewcommand{\1}{\mathbbm{1}}
\newcommand{\e}{\mbox{\rm e}}

\makeatother
\begin{document}

\title[A multifractal analysis for Stern-Brocot intervals]{A
multifractal analysis for Stern-Brocot intervals, continued fractions
and Diophantine growth rates}
\author{Marc Kesseb\"ohmer}
\address{Fachbereich 3 -- Mathematik und Informatik, Universit\"at Bremen,
Bibliothekstr. 1, D--28359 Bremen, Germany}
\email{mhk@math.uni-bremen.de}

\author{Bernd O. Stratmann} \address{Mathematical Institute, University
of St. Andrews, North Haugh, St. Andrews KY16 9SS, Scotland}
\email{bos@maths.st-and.ac.uk}


\date{}


\begin{abstract}
In this paper we obtain multifractal generalizations of classical
results by L{\'e}vy and Khintchin in metrical Diophantine approximations
and measure theory of continued fractions. We give a complete
multifractal
analysis for Stern--Brocot intervals, for continued fractions and
for certain Diophantine growth rates. In particular, we give detailed
discussions of two multifractal spectra
closely related to the Farey map and to the Gauss map.
\end{abstract}

\maketitle

\section{Introduction and statements of result}

In this paper we give a multifractal analysis for Stern--Brocot
intervals,
continued fractions and certain Diophantine growth rates. We apply
and extend the multifractal formalism for average growth rates of
\cite{KesseboehmerStratmann:04a} to obtain a complete multifractal
description
   of two dynamical systems originating from the set of real numbers.

Recall that the process of writing an element $x$ of the unit interval
in its regular continued fraction expansion \[
x=[a_{1}(x),a_{2}(x),a_{3}(x),\ldots]=\cfrac{1}{a_{1}(x)+\cfrac{1}{a_{2}(x)+\cfrac{1}{a_{3}(x)+\cdots}}}\]
can be represented either by a uniformly hyperbolic dynamical system
which is based on an infinite alphabet and hence has infinite
topological
entropy, or by a non-uniformly hyperbolic dynamical system based on
a finite alphabet and having finite topological entropy. Obviously,
for these two systems the standard theory of multifractals (see e.g.
\cite{Pesin:97}) does not apply, and therefore it is an interesting task
to give a multifractal analysis for these two number theoretical
dynamical systems. There is a well known result which gives some
information
in the generic situation, that is for a set of full $1$-dimensional
Lebesgue measure $\lambda$. Namely with
$p_{n}(x)/q_{n}(x):=[a_{1}(x),a_{2}(x),\ldots,a_{n}(x)]$
referring to the $n$-th approximant of $x$, we have for $\lambda$-almost
every $x\in[0,1)$, \[
\ell_{1}(x):=\lim_{n\rightarrow\infty}\frac{2\log
q_{n}(x)}{\sum_{i=1}^{n}a_{i}(x)}=0.\]
   Note that by employing the analogy between regular continued fraction
expansions of real numbers and geodesics on the modular surface, the
number $2\log q_{n}(x)$ can be interpreted as the 'hyperbolic length'
associated with the approximant $p_{n}(x)/q_{n}(x)$. Also,
the parameter $n$ represents the word length associated with
$p_{n}(x)/q_{n}(x)$
with respect to the dynamical system on the infinite alphabet, whereas
$\sum_{i=1}^{n}a_{i}(x)$ can be interpreted as the word length
associated
with $p_{n}(x)/q_{n}(x)$ with respect to the dynamical system on
the finite alphabet.  There are two classical results by Khintchin and
L{\'e}vy \cite{Levy:29}, \cite{Levy:36}, \cite{Khinchin:35}, \cite{Khihtchine:36} which allow a
closer inspection of the limit $\ell_{1}$.  That is, for
$\lambda$-almost every $x\in[0,1)$ we have, with
$\chi:=\pi^{2}/(6\log2)$,
\[
\ell_{2}(x):=\lim_{n\rightarrow\infty}\frac{\sum_{i=1}^{n}a_{i}(x)}{n}=\infty
 \textrm{  and } \ell_{3}(x):=\lim_{n\rightarrow\infty}\frac{2\log
q_{n}(x)}{n}=\chi.\]
   Clearly, dividing the sequence in $\ell_{3}$ by the sequence in
$\ell_{2}$ leads to the sequence in $\ell_{1}$. Therefore, if we
define the level sets \[
\mathcal{L}_{i}(s):=\left\{ x\in[0,1):\ell_{i}(x)=s\right\}, 
\textrm{  for  }  s\in\R\,,\]
   then these classical results by L{\'e}vy and Khintchin imply for the
Hausdorff dimensions ($\dim_{H}$) of these level sets \[
\dim_{H}(\mathcal{L}_{1}(0))=\dim_{H}(\mathcal{L}_{2}(\infty)\cap\mathcal{L}_{3}(\chi))=1.\]
A natural question to ask is what happens to this relation between
these Hausdorff dimensions for prescribed non-generic limit behavior.
Our first main results in this paper will give an answer to this
question.
Namely, with $\gamma:=(1+\sqrt{5})/2$ referring to the Golden
Mean, we show that for each $\alpha\in[0,2\log\gamma]$ there exists
a number $\alpha^{\sharp}=\alpha^{\sharp}(\alpha)\in\R\cup\{\infty\}$
such that, with the convention $\alpha^{\sharp}(0):=\infty$ and
$0\cdot\alpha^{\sharp}(0):=\chi$,
\[
\dim_{H}(\mathcal{L}_{1}(\alpha))=\dim_{H}(\mathcal{L}_{2}(\alpha^{\sharp})\cap\mathcal{L}_{3}
(\alpha\cdot\alpha^{\sharp})).\]
Furthermore, for the dimension function $\tau$ given by \[
\tau(\alpha):=\dim_{H}(\mathcal{L}_{1}(\alpha)),\]
   we show that $\tau$ can be expressed explicitly in terms of the
Legendre transform $\widehat{P}$ of a certain pressure function
$P$, referred to as the Stern--Brocot pressure. For the function $P$
we obtain the result that it is real-analytic on the interval
$(-\infty,1)$
and vanishes on the complement of this interval. We then show that
the dimension function $\tau$ is continuous and strictly decreasing
on $[0,2\log\gamma]$, that it vanishes outside the interval
$[0,2\log\gamma)$,
and that for $\alpha\in[0,2\log\gamma]$ we have \[
\alpha\cdot\tau(\alpha)=-\widehat{P}(-\alpha).\]
Before we state the main theorems, let us recall the following classical
construction of Stern--Brocot intervals (cf. \cite{Stern1858},
\cite{Brocot:1860}).
For each $n\in\N_{0}$, the elements 
of the $n$-th member of the Stern--Brocot sequence \[
\mathfrak{T}_{n}:=\left\{
\frac{s_{n,k}}{t_{n,k}} : k=1,\ldots,2^{n}+1\right\} \]
   are defined recursively as follows.
\begin{itemize}
\item $s_{0,1}:=0\,\,$ and $\,\, s_{0,2}:=t_{0,1}:=t_{0,2}:=1$;
\item $s_{n+1,2k-1}:=s_{n,k}\quad\textrm{and}\quad
t_{n+1,2k-1}:=t_{n,k},$
for $k=1,\ldots,2^{n}+1$;
\item $s_{n+1,2k}:=s_{n,k}+s_{n,k+1}\quad\textrm{and}\quad
t_{n+1,2k}:=t_{n,k}+t_{n,k+1}$,
for $k=1,\ldots2^{n}$.
\end{itemize}
With this ordering of the rationals in $[0,1]$ we define the set
$\mathcal{T}_{n}$ of Stern--Brocot intervals of order $n$ by \[
\mathcal{T}_{n}:=\left\{
T_{n,k}:=\left[\frac{s_{n,k}}{t_{n,k}},\frac{s_{n,k+1}}{t_{n,k+1}}\right):\,
k=1,\ldots,2^{n}\right\} .\]
Clearly, for each $n\in\N_{0}$ we have that $\mathcal{T}_{n}$ represents
a partition of the interval $[0,1)$. The first members in this sequence
of sets are the following, and it should be clear how to continue
this list using the well known method of mediants.
\begin{eqnarray*}
\mathcal{T}_{0}= & \left\{ \left[\frac{0}{1},\frac{1}{1}\right)\right\}
\\
\mathcal{T}_{1}= & \left\{
\left[\frac{0}{1},\frac{1}{2}\right),\left[\frac{1}{2},\frac{1}{1}\right)\right\}
\\
\mathcal{T}_{2}= & \left\{
\left[\frac{0}{1},\frac{1}{3}\right),\left[\frac{1}{3},\frac{1}{2}\right),
\left[\frac{1}{2},\frac{2}{3}\right),\left[\frac{2}{3},\frac{1}{1}\right)\right\}
\\
\mathcal{T}_{3}= & \left\{
\left[\frac{0}{1},\frac{1}{4}\right),\left[\frac{1}{4},\frac{1}{3}\right),
\left[\frac{1}{3},\frac{2}{5}\right),\left[\frac{2}{5},\frac{1}{2}\right),
\left[\frac{1}{2},\frac{3}{5}\right),\left[\frac{3}{5},\frac{2}{3}\right),
\left[\frac{2}{3},\frac{3}{4}\right),\left[\frac{3}{4},\frac{1}{1}\right)\right\}
\\
\vdots & \vdots\end{eqnarray*}
As already mentioned above,  our multifractal analysis will make use
of the Stern--Brocot pressure function $P$. This function is defined for
$\theta\in\R$
by \[
P(\theta):=\lim_{n\rightarrow\infty}\frac{1}{n}\log\sum_{T\in\mathcal{T}_{n}}
\left|T\right|^{\theta}.\]
In here, $\left|T\right|$ refers to the Euclidean length of the interval $T$.
We will see that  $P$ is a well--defined convex function (cf. Proposition
\ref{pro:AnalyticPropertiesP}). Also, note that we
immediately have that\[
P(\theta)=\lim_{n\rightarrow\infty}\frac{1}{n}\log\sum_{k=1}^{2^{n}}
\left(\frac{s_{n,k+1}}{t_{n,k+1}}-
\frac{s_{n,k}}{t_{n,k}}\right)^{\theta}=\lim_{n\rightarrow\infty}
\frac{1}{n}\log\sum_{k=1}^{2^{n}}\left(\frac{1}{t_{n,k}\cdot
t_{n,k+1}}\right)^{\theta}.\]
The following theorem gives the first main results of this paper.
In here, $\widehat{P}$ refers to the Legendre transform of $P$,
given for $t \in\R$ by
$\widehat{P}(t):=\sup_{\theta\in\R}\{\theta t-P(\theta)\}$.

\begin{figure}
\psfrag{alpha-}{$-2\theta\log(\gamma)$}
\psfrag{DimFG}{$\dim_H({\mathcal{L}}_1(\alpha))$}
\psfrag{pss}{$P(\theta)$} \psfrag{2 log g}{$2\log(\gamma)$}
\psfrag{1}{\(1\)} \psfrag{a}{$\alpha$} \psfrag{s}{$\theta$}
\includegraphics[%
    width=1.0\columnwidth,
    keepaspectratio]{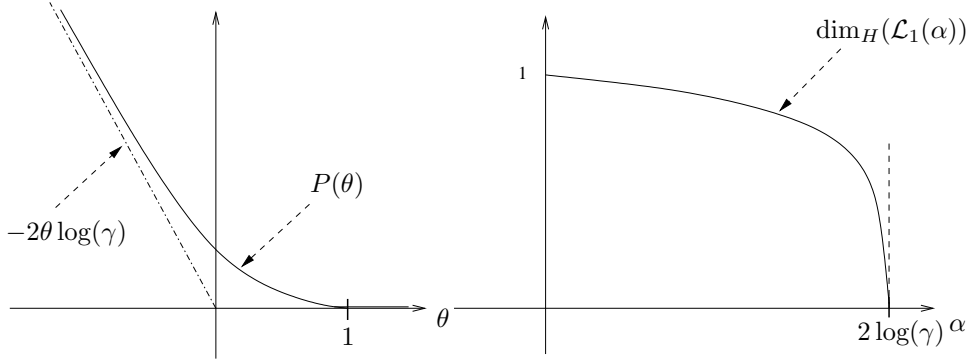}

\caption{The Stern--Brocot pressure $P$ and the multifractal
spectrum $\tau$ for $\ell_1$.\label{cap:The-Stern--Brocot-pressure}}
\end{figure}

\begin{thm}
\label{Thm:main}
{\em (see Fig. \ref{cap:The-Stern--Brocot-pressure})} 
\begin{enumerate}
\item The Stern--Brocot pressure $P$ is convex, non-increasing and
differentiable
throughout $\R$. Furthermore, $P$ is real--analytic on the interval
$(-\infty,1)$ and is equal to $0$ on $[1,\infty)$.
\item For every $\alpha\in[0,2\log\gamma]$, there exist
$\alpha^{*}=\alpha^{*}(\alpha)\in\R$
and $\alpha^{\sharp}=\alpha^{\sharp}(\alpha)\in\R\cup\{\infty\}$
related by $\alpha\cdot\alpha^{\sharp}=\alpha^{*}$ such that, with
the conventions $\alpha^{*}(0):=\chi$ and $\alpha^{\sharp}(0):=\infty$,
\begin{eqnarray*}
\dim_{H}\left(\mathcal{L}_{1}(\alpha)\right) & = &
\dim_{H}\left(\mathcal{L}_{2}(\alpha^{\sharp})\cap\mathcal{L}_{3}(\alpha^{*})\right)
\,\,\left(=:\tau(\alpha)\right).\end{eqnarray*}
Furthermore, the dimension function $\tau$ is continuous and strictly
decreasing on $[0,2\log\gamma]$, it vanishes outside the interval
$[0,2\log\gamma)$, and for $\alpha\in[0,2\log\gamma]$ we have\[
\alpha\cdot\tau(\alpha)=-\widehat{P}(-\alpha),\]
   where $\tau(0):=\lim_{\alpha\searrow0}-\widehat{P}(-\alpha)/\alpha=1$.
Also, for the left derivative of $\tau$ at $2\log\gamma$ we have
$\lim_{\alpha\nearrow2\log\gamma}\tau'\left(\alpha\right)=-\infty$.
\end{enumerate}
\end{thm}
Theorem \ref{Thm:main} has some interesting implications for other 
canonical level sets. In order to 
state these, recall that the elements
of $\mathcal{T}_{n}$ cover the interval $[0,1)$ without overlap.
Therefore, for each $x\in[0,1)$ and $n\in\N$ there exists a unique
Stern--Brocot interval $T_{n}(x)\in\mathcal{T}_{n}$ containing $x$.
The interval $T_{n}(x)$ is covered by two neighbouring intervals from
$\mathcal{T}_{n+1}$, a left and a right subinterval.  If $T_{n+1}(x)$
is the left of these then we encode this event by the letter $A$,
otherwise we encode it by the letter $B$.  In this way every
$x\in[0,1)$ can be described by a unique sequence of nested
Stern--Brocot intervals of any order that contain $x$, and therefore
by a unique infinite word in the alphabet $\{ A,B\}$.  It is well
known that this type of coding is canonically associated with the
continued fraction expansion of $x$ (see Section 2
 for the details).  In particular,
this allows to relate the level sets $\mathcal{L}_{1}$ and
$\mathcal{L}_{3}$ to level sets given by means of the Stern--Brocot
growth rate $\ell_{4}$ of the nested sequences
$\left(T_{n}(x)\right)$, and to level sets of certain Diophantine
growth rates $\ell_{5}$ and $\ell_{6}$.  These growth
rates are given by (assuming that these limits exist) \[
\ell_{4}(x):=\lim_{n\rightarrow\infty}\frac{\log\left|T_{n}(x)\right|}{-n},\,\,\]
\[
\ell_{5}(x):=\lim_{n\rightarrow\infty}\frac{2\log\left|x-\frac{p_{n}(x)}{q_{n}(x)}\right|}
{-\sum_{i=1}^{n}a_{i}(x)}\quad\textrm{and}\quad\ell_{6}(x):={\displaystyle
\lim_{n\rightarrow\infty}\frac{2\log\left|x-\frac{p_{n}(x)}{q_{n}(x)}\right|}{-n}.}\]

\begin{figure}
\psfrag{alpha-}{$-2\theta\log\gamma$}\psfrag{2loggam}{$2\log\gamma$}
\psfrag{DimFG}{$\dim_H({\mathcal{L}}_3(\alpha))$} \psfrag{pss}{$P_D
(\theta)$}
\psfrag{1}{\(1\)}\psfrag{1}{\(1\)}\psfrag{1/2}{\(\frac{1}{2}\)}\psfrag{a}{$\alpha$}
\psfrag{s}{$\theta$} \psfrag{chi}{$\chi:=\frac{\pi^2}{6\log2}$} \psfrag{d0}{$d_0$}

\includegraphics[%
    width=1.0\columnwidth,
    keepaspectratio]{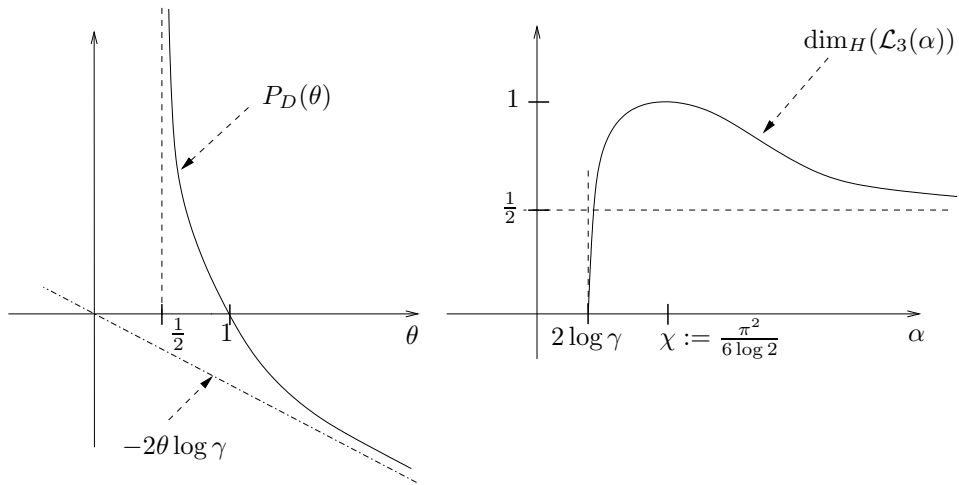}

\caption{The Diophantine pressure $P_{D}$ and the
multifractal
spectrum $\tau_{D}$ for $\ell_{3}$.\label{cap:Diaophantine-pressure}}
\end{figure}

\begin{prop}
\label{Thm:2} Let $x \in [0,1)$ be given. If one of the limits in 
$\{ \ell_{1}(x), \ell_{4}(x), \ell_{5}(x)\}$ exists then also the 
other two do exist, and   \[
\ell_{1}(x)=\ell_{4}(x)=\ell_{5}(x).\]
Furthermore, $\ell_{3}(x)$ exists if and only if $\ell_{6}(x)$ exists, 
and if one of these exists then
\[ \ell_{3}(x)=\ell_{6}(x).\]
   By Theorem \ref{Thm:main}, it therefore follows that for each
$\alpha\in[0,2\log\gamma]$,
\[
\dim_{H}\left(\mathcal{L}_{4}(\alpha)\right)=\dim_{H}\left(\mathcal{L}_{5}(\alpha)\right)=
\dim_{H}\left(\mathcal{L}_{2}(\alpha^{\sharp})\cap\mathcal{L}_{6}(\alpha^{*})\right)=\tau(\alpha).\]
\end{prop}
Note that the level sets $\mathcal{L}_{4}(\alpha)$ have already been 
under consideration in   
\cite{KesseboehmerStratmann:04}. There they were introduced in terms 
of homological growth rates of hyperbolic geodesics (see Remark 
\ref{rem:Vergleich} (2)).  Clearly, Theorem \ref{Thm:main} and Proposition \ref{Thm:2} consider the
dynamical system associated with the finite alphabet, a system  which is closely
related to the Farey map. Now, our second main result gives a multifractal
analysis for the system based on the infinite alphabet, and this system is
closely related
to the Gauss map. In here, the relevant pressure function is the
\emph{Diophantine
pressure} $P_{D}$, which is given by \[
P_{D}(\theta):=\lim_{k\rightarrow\infty}\frac{1}{k}\log\sum_{\left[a_{1},\ldots,a_{k}\right]}
q_{k}\left(\left[a_{1},\ldots,a_{k}\right]\right)^{-2\theta}, \textrm{  
for  } \theta>\frac{1}{2}.\]
We remark that a very detailed analysis of the function $P_{D}$ can be 
found in \cite{Mayer}. Our second main result is the following. 
\begin{thm}
\label{thm:main3} {\em (see Fig. \ref{cap:Diaophantine-pressure})} The function $P_{D}$ has a singularity
at $1/2$, and
$P_{D}$ is decreasing, convex and real-analytic on $\left(1/2,\infty\right)$.
Furthermore, for $\alpha\in[2\log\gamma,\infty)$ we have\[
\dim_{H}\left(\mathcal{L}_{3}(\alpha)\right)=\dim_{H}\left(\mathcal{L}_{6}(\alpha)\right)=
\frac{\widehat{P}_{D}(-\alpha)}{-\alpha}=:\tau_{D}(\alpha).\]
Also, the dimension function $\tau_{D}$ is real-analytic on
$(2\log\gamma,\infty)$,
it is increasing on $[2\log\gamma,\chi]$ and decreasing on
$[\chi,\infty)$.
In particular, $\tau_{D}$ has a point of inflexion at some point greater
than
$\chi$ and a unique maximum equal to $1$ at $\chi$. Additionally,
$\lim_{\alpha\to\infty}\tau_{D}\left(\alpha\right)=1/2$,
$\lim_{\alpha\searrow2\log\gamma}\tau_{D}\left(\alpha\right)=0$
and $\lim_{\alpha\searrow2\log\gamma}\tau_{D}'(\alpha)=\infty$.
\end{thm}

The paper is organized as follows. In Section 2 we first recall two
ways of coding elements of the unit interval. One is based on a finite 
alphabet and
the other on an infinite alphabet, and both are defined in terms of  the modular group.
These codings are
canonically related to regular continued fraction expansions, and
we end the section by commenting on a 1-1 correspondence between
Stern--Brocot
sequences and finite continued fraction expansions. In Section 3 we
introduce certain cocycles which are relevant in our multifractal
analysis. In particular, we give various estimates relating these
cocycles with the geometry of the modular codings and with the sizes
of the Stern--Brocot intervals. This will then enable us to prove the
first part of Proposition \ref{Thm:2}. Section 4 is devoted to the
discussion
of several aspects of the Stern--Brocot pressure and its Legendre
transform.
In Section 5 we give the proof of Theorem \ref{Thm:main}, which we
have split into the parts \emph{The lower bound}, \emph{The upper 
bound},
and \emph{Discussion of boundary points of the spectrum}.  Finally, in Section 6 we
give the proof of Theorem \ref{thm:main3} by showing how to adapt
the multifractal formalism developed in Section 4 and 5 to the 
situation here.  

Throughout, we shall use the notation $f\ll g$ to denote that
for two non-negative functions $f$ and $g$ we have that $f/g$
is uniformly bounded away from infinity. If $f\ll g$ and $g\ll f$,
then we write $f\asymp g$.

\begin{rem}
 One immediately verifies that the results of Theorem
  \ref{Thm:main} and
Proposition \ref{Thm:2} can be expressed in terms of the Farey map $\mathfrak{f}$ 
acting on $[0,1]$, and then 
$\tau$ represents the multifractal spectrum of
the measure  of maximal entropy (see e.g. \cite{Nakaishi:00}).
Likewise, the results of Theorem \ref{thm:main3} can
be written in terms of the Gauss map $\mathfrak{g}$, and then in
this terminology $\tau_D$ describes the Lyapunov spectrum of $\mathfrak{g}$. 
For the definitions of $\mathfrak{f}$  and $\mathfrak{g}$ and for a discussion
of their relationship we refer to Remark \ref{rem:FG}.

\end{rem}
\begin{rem}
Since the theory of multifractals started through essays of Mandelbrot
\cite{Mandelbrot:74} \cite{Mandelbrot:88}, Frisch and Parisi
\cite{FrischParisi:85},
and Halsey et al. \cite{Halsey:86}, there has been a steady increase of
the literature on multifractals and calculations of specific
multifractal
spectra. For a comprehensive account on the mathematical work we refer
to  \cite{PesinWeiss:96} and \cite{Pesin:97}. Essays which are
   closely related to the
work on multifractal number
theory in this paper are for instance \cite{Byrne:98},  
\cite{FengOlivier:03}, \cite{KesseboehmerStratmann:04},
\cite{HanusMauldiUrban:02}, \cite{MesonVericat:04}, \cite{Nakaishi:00}
and \cite{PW}. We remark that brief sketches of some parts of Theorem  \ref{thm:main3}
have already been given in \cite{KesseboehmerStratmann:04}. The
results there do for 
instance not cover the boundary points of the spectra. 
Furthermore,  note that for the $\ell_{6}$--spectrum  partial results have been 
established in \cite{PW} (Corollary 2).
\end{rem}

\section{The Geometry of Modular Codings by Finite and Infinite
Alphabets\label{sec:coding}}
Let $\Gamma:=\textrm{PSL}_{2}\left(\Z\right)$ refer to the modular
group acting on the upper half-plane $\H$. It is well--known that
$\Gamma$ is generated by the two elements $P$ and $Q$, given by\[
P:z\mapsto z-1\,\,\,\textrm{and}\,\,\, Q:z\mapsto\frac{-1}{z}.\]
\begin{figure}[ht]
\psfrag{z0''}{$z_0'$} \psfrag{F}{$F$}\psfrag{R(F)}{$R(F)$}
\psfrag{R2(F)}{$R^2(F)$}\psfrag{i}{$i$} \psfrag{0}{$0$}\psfrag{1}{$1$}

\includegraphics[%
    width=0.50\columnwidth,
    keepaspectratio]{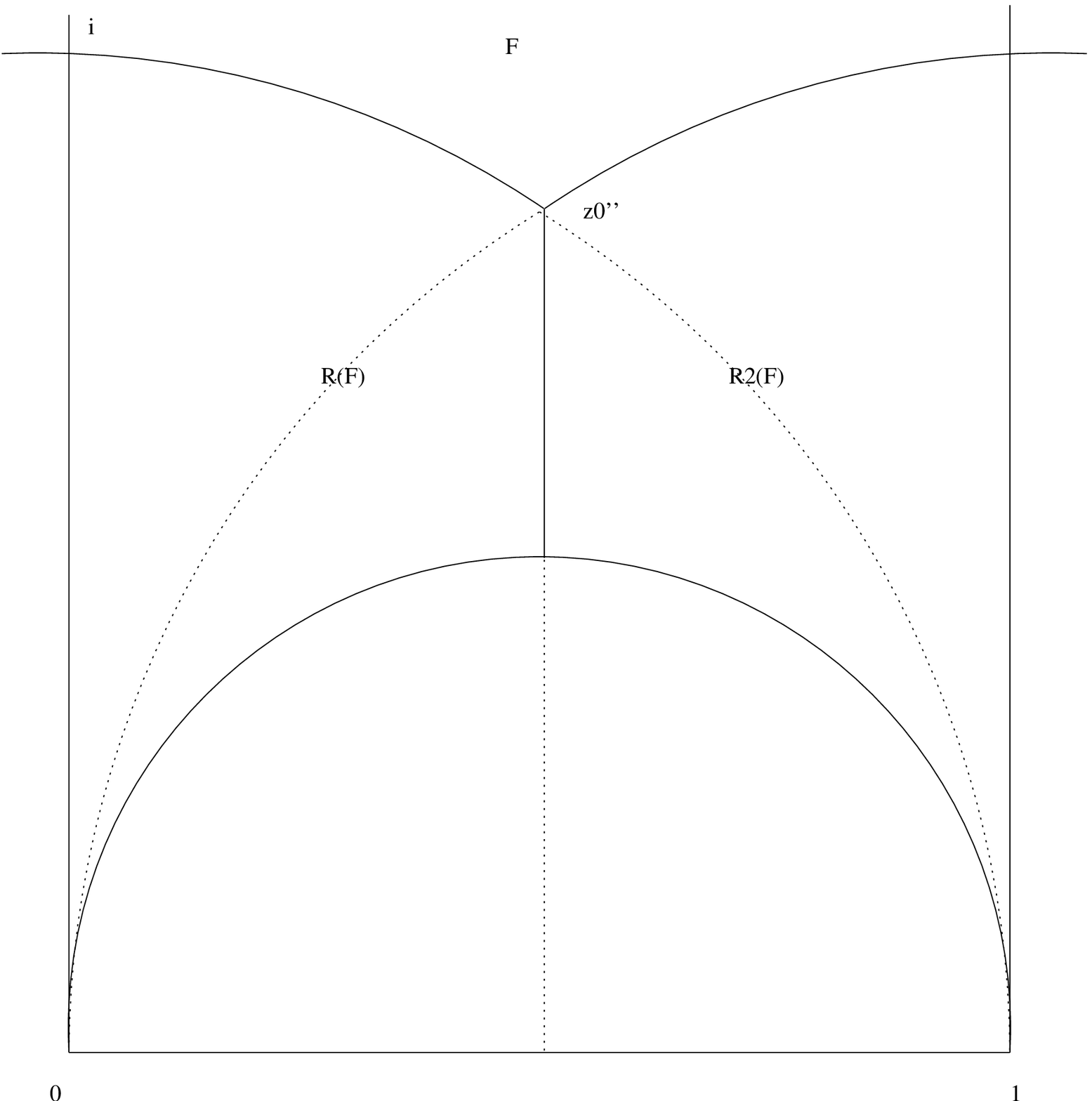}

\caption{A fundamental domain $F$ for
$\textrm{PSL}{}_{2}\left(\Z\right)$
and the images under $R$ and $R^{2}$. \label{cap:FundDomain}}
\end{figure}

\begin{figure}[ht]
\psfrag{z0''}{ $z_0'$}\psfrag{z0}{ $z_0$} \psfrag{A(z0)} {$A(z_0)$}
\psfrag{AA(z0)}{ $AA(z_0)$} \psfrag{AAA(z0)}{ $AAA(z_0)$}
\psfrag{BB(z0)}{$BB(z_0)$}\psfrag{BBB(z0)}{$BBB(z_0)$}
\psfrag{AB(z0)}{\tiny $AB(z_0)$} \psfrag{BA(z0)}{\tiny $BA(z_0)$}
\psfrag{B(z0)}{ $B(z_0)$} \psfrag{1/2}{$1/2$}
\psfrag{1/3}{$1/3$} \psfrag{1/4}{$1/4$} \psfrag{2/3}{$2/3$}
\psfrag{1/4}{$1/4$} \psfrag{3/4}{$3/4$}
\psfrag{3/5}{$3/5$}\psfrag{2/5}{$2/5$} \psfrag{0}{$0$}
\psfrag{1}{$1$}
\psfrag{T}{$\times$}\psfrag{T22}{$T_{2,2}$}\psfrag{T21}{$T_{2,1}$}\psfrag{T23}{$T_{2,3}$}\psfrag{T24}{$T_{2,4}$}

\psfrag{T32}{$T_{3,2}$}\psfrag{T31}{$T_{3,1}$}\psfrag{T33}{$T_{3,3}$}\psfrag{T34}{$T_{3,4}$}

\psfrag{T35}{$T_{3,5}$}\psfrag{T36}{$T_{3,6}$}\psfrag{T37}{$T_{3,7}$}\psfrag{T38}{$T_{3,8}$}

\psfrag{F0}{$F_0$}\psfrag{ABF0}{$A(F_0)=B(F_0)$}

\includegraphics[%
    width=0.98\columnwidth,
    keepaspectratio]{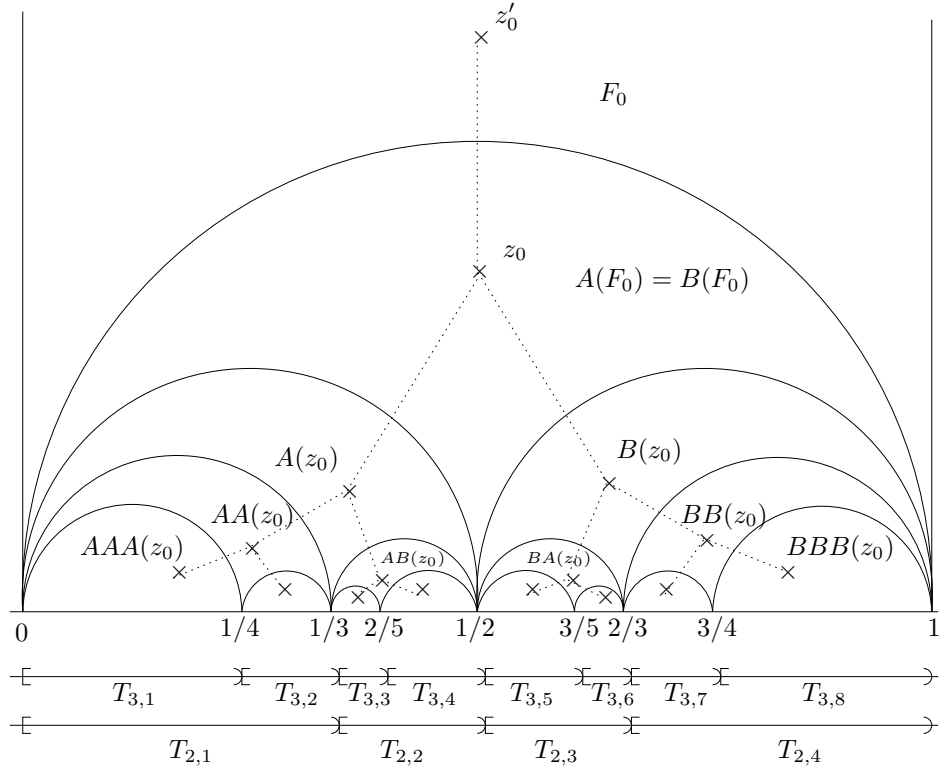}

\caption{Part of the Cayley graph rooted at $z_{0}$, for
$\Gamma_0(z_0)$
restricted to $[0,1]\times\R^{+}$, and the
Stern--Brocot intervals of order $2$ and $3$. \label{cap:fig1}}
\end{figure}

Defining relations for $\Gamma$ are $Q^{2}=(PQ)^{3}=\{\textrm{id}.\}$,
and a fundamental domain $F$ for $\Gamma$ is the hyperbolic
quadrilateral
with vertices at $i,1+i,\{\infty\}$ and $z_{0}':=(1+i\sqrt{3})/2$.
For $R:=QP$ such that $R:z\mapsto-1/(z-1)$, one easily verifies
that $\Gamma_{0}:=\Gamma/\left\langle R\right\rangle $ is a subgroup
of $\Gamma$ of index $3$ and that $F_{0}$ is a fundamental domain
for $\Gamma_{0}$, for $F_{0}:=F\cup R(F)\cup R^{2}(F)$ the ideal
triangle with vertices at $0,1$ and $\{\infty\}$ (see Fig.
\ref{cap:FundDomain}).
Consider the two elements $A,B\in\Gamma$ given by\[
A:=\left(Q^{-1}PQ\right):z\mapsto\frac{z}{z+1}\quad\textrm{and}\quad
B:=\left(P^{-1}A^{-1}P\right):z\mapsto\frac{-1}{z-2},\]
   and let $G$ denote the free semi-group generated by $A$ and $B$.
It is easy to see that for $z_{0}:=A(z_{0}')=B(z_{0}')=(1+i/\sqrt{3})/2$
we have that the Cayley graph of $G$ with respect to $z_{0}$ coincides
with the restriction to  \[\left\{
z\in\H{:}\,0\leq\Re(z)\leq1,\;0<\Im(z)\leq1/2\right\} \]
of the the Cayley graph of $\Gamma_{0}$ with respect to $z_{0}$
(see Fig. \ref{cap:fig1}).

\subsection{Finite Coding} 

Let $\Sigma:=\left\{ A,B\right\} ^{\N}$ denote the full shift space
on the finite alphabet $\left\{ A,B\right\}$, for $A,B \in \Gamma$ 
given as above. Also, let
$\Sigma$ be equipped with the usual left-shift
$\sigma:\Sigma\rightarrow\Sigma$.
Then $\Sigma$ is clearly isomorphic to the completion of
$G$, where the completion is taken with respect to a suitable metric
on $G$ (see \cite{Floyd:80}). One immediately verifies that the
canonical map \begin{eqnarray*}
\pi:\qquad\Sigma\qquad & \rightarrow & \qquad[0,1],\\
(x_{1},x_{2},\ldots) & \mapsto & \lim_{n\rightarrow\infty}x_{1}\cdots
x_{n}(z_{0}),\end{eqnarray*}
is 1--1 almost everywhere, in the sense that it is  2--1 on the rationals in $[0,1]$ and 1--1 on $\I$,
for $\I$ referring to the irrational numbers in $\left[0,1\right]$. Note that 
the Stern--Brocot sequence $\mathfrak{T}_{n+1}$ coincides
with the set of vertices at infinity of $\{ g(F_{0}):g\in
G$ of word length $n\} $, for each
$n\in\N$.

\subsection{Infinite Coding}\label{2.2}

For the infinite alphabet $\left\{
X^{n}: n\in\N, X \in \{A,B\}\right\} $ we define the shift space of
finite type \[ \Sigma^{*}:=\left\{ \left(X^{n_{1}},
Y^{n_{2}},X^{n_{3}},\ldots\right):\{X,Y\}=\{A,B\}, n_i \in \N \hbox{ 
for all } i \in \N\right\}, \] which we assume to be equipped with the usual
left-shift $\sigma^{*}:\Sigma^{*}\rightarrow\Sigma^{*}$.  Then there
exists a canonical bijection $\pi^{*}$, given by\begin{eqnarray*}
\pi^{*}:\qquad\Sigma^{*}\qquad & \rightarrow & \I\\
(y_{1},y_{2},\ldots) & \mapsto &
\lim_{k\rightarrow\infty}y_{1}y_{2}\cdots y_{k}(z_{0}).\end{eqnarray*}
   This coding is closely related to the continuous fraction expansion.
   Namely, if
$y=\left(X^{n_{1}},Y^{n_{2}},X^{n_{3}},\ldots\right)$ then \[
\pi^{*}\left(y\right)=\left\{ \begin{array}{lll}
[n_{1}+1,n_{2},n_{3},\ldots] & \,\,\textrm{for} &
X=A\\
{}[1,n_{1},n_{2},\ldots] & \,\,\textrm{for} &
X=B.\end{array}\right.\]
Also, if $S:[0,1]\rightarrow[0,1]$ and
$s:\Sigma^{*}\rightarrow\Sigma^{*}$
are given by, for  $x \in [0,1]$ and $\left\{ X,Y\right\} =\left\{
A,B\right\}$,
\[ S(x) :=\left(1-x\right) \quad\textrm{and}\quad
s\left(X^{n_{1}},
Y^{n_{2}},X^{n_{3}},\ldots\right) :=\left(Y^{n_{1}},
X^{n_{2}},Y^{n_{3}},\ldots\right),\]
then by symmetry we have that $S\circ\pi^{*}=\pi^{*}\circ s$.

To overcome the fact that $\left(\Sigma^{*},\sigma^{*}\right)$ is
not topological transitive, we also require the full shift space
$\left(\overline{\Sigma},\overline{\sigma}\right)$ over $\N$. In here, 
$\overline{\Sigma} := \N^{\N}$ and $\overline{\sigma}$ refers to
the left-shift map on $\overline{\Sigma} $. Clearly,
$\left(\overline{\Sigma},\overline{\sigma}\right)$
is \emph{finitely primitive} in the sense of \cite{MauldinUrbanski:03},
and we remark that this property is a necessary preliminary for
the thermodynamical formalism
 used throughout this paper.

Note that the two shift spaces $(\Sigma^{*}, \sigma^{*})$ and 
$\left(\overline{\Sigma},\overline{\sigma}\right)$ are related by the  2-1 factor 
map $p$, which is given by
\begin{equation}
p:\Sigma^{*}\to\overline{\Sigma},\,
\,\,\,\left(X^{n_{1}},Y^{n_{2}},X^{n_{3}},\ldots\right)\mapsto\left(n_{1},n_{2},
n_{3},\ldots\right).\label{eq:FactorMap}\end{equation}

\begin{rem}\label{rem:FG}
Note that the finite coding is in 1-1 correspondence to the coding
of $[0,1]$ via the inverse branches $f_{1}$ and $f_{2}$ of the \emph{Farey
map} $\mathfrak{f}$. In here, $f_{1}$ and $f_{2}$ are given by $f_{1}(x)=x/(x+1)$ and
$f_{2}(x)= 1/(x+1)$, for $x \in [0,1]$.
One easily verifies that $f_{1}=A$ and $f_{2} \circ S=B$, and hence
$\Sigma$ can be interpreted as arising from a `twisted Farey
map'. Similarly, one notices that $\Sigma^{*}$ is closely related
to
the coding of $[0,1]$ via the infinitely many branches of the Gauss
map $\mathfrak{g}$, which is given by $\mathfrak{g}(x):=1/x \;\mathrm{mod }\;1$
for $x \in [0,1]$. More precisely, we have that the dynamical
system $(\I,\mathfrak{g})$  is a topological 2--1 factor of the dynamical system
 $(\Sigma^*,\sigma^*)$, where the factor map can be established either on the symbolic
 level via $p$ or on the geometric level via ${\mathfrak{f}}$. 
The situation is summarized  in the  following commuting diagram. 
{\Large \begin{eqnarray*}\xymatrix{
\I \ar[d]_{\mathfrak{g}} && \ar[ll]_{\pi_{\mathrm{CF}}}\overline{\Sigma}\ar[d]_{\overline{\sigma}}&&
\ar[ll]_{p}\Sigma^* \ar[rr]^{\pi^*}\ar[d]_{\sigma^*} &&
 \I \ar[d]^{{\mathfrak{g}_s}} \ar[rr]^{\mathfrak{f}} && \I \ar[d]^{{\mathfrak{g}}} \\
\I && \ar[ll]^{\pi_{\mathrm{CF}}} \overline{\Sigma} &&
\ar[ll]^{p}\Sigma^* \ar[rr]_{\pi^*}  && \I   \ar[rr]_{\mathfrak{f}} && \I \\}  
\end{eqnarray*}}
In here, $\mathfrak{g}_s := \pi^*\circ\sigma^*\circ(\pi^*)^{-1}$ denotes the 'twisted Gauss map' and 
$\pi_{\mathrm{CF}}$ is given by
$\pi_{\mathrm{CF}}(n_1,n_2,\ldots ):=[n_1,n_2,\ldots ]$  for $(n_1,n_2,\ldots )\in\overline\Sigma$.
 Note that both $\pi_{\mathrm{CF}}$ and $\pi^*$ are bijections and that $\pi_{\mathrm{CF}}\circ p=\mathfrak{f}\circ\pi^*$. 
\end{rem}
\subsection{Stern--Brocot sequences versus continued fractions}

We end this section by showing that there is a 1--1 correspondence
between
the elements of the Stern--Brocot sequence and finite continued fraction
expansions. This will turn out to be useful in the sequel.

For $n\geq2$, let $A_{k}^{n}$ refer to the set all $k$-tuples of
positive integers which add up to $n$ and whose $k$-th entry exceeds
$1$.
That is, \begin{equation}
A_{k}^{n}:=\left\{ \left(a_{1},a_{2},\ldots,
a_{k}\right)\in\mathbb{N}^{k}:\sum_{i=1}^{k}a_{i}=n,\,\,
a_{k}\not=1\right\} .\label{eq:DefAnk}\end{equation}
   Since $a_{k}\not\neq1$, we can identify an element $\left(a_{1},
\ldots,a_{k}\right)\in A_{k}^{n}$
in a unique way with the finite continued fraction expansion
$\left[a_{1},a_{2},\ldots a_{k}\right]$.
Also, one easily verifies that for $1\leq k\leq n-1$,\begin{equation}
\card\left(A_{k}^{n}\right)={\binom{n-2}{
k-1}}.\label{eq:CardAkn}\end{equation}

\begin{lem}
\label{lem:TnbyCF} For all $n\geq2$ we have\[
\bigcup_{k=1}^{n-1}\bigcup_{(a_{1},a_{2},\ldots,a_{k}) \in A_{k}^{n}}\left[a_{1},a_{2},\ldots a_{k}
\right]=\mathfrak{T}_{n-1}\setminus\mathfrak{T}_{n-2}=
\left\{ \frac{s_{n-1,2\ell}}{t_{n-1,2\ell}}:1\leq\ell\leq2^{n-2}\right\}
.\]
Furthermore, if $\left(s_{n,k}/t_{n,k}\right)=\left[a_{1},a_{2},
\ldots,a_{m}\right]\in\mathfrak{T}_{n}\setminus\mathfrak{T}_{n-1}$ then its two siblings
in $\mathfrak{T}_{n+1}\setminus\mathfrak{T}_{n}$ are, for $\left\{ u,v\right\} =\left\{
2k,2k-2\right\} $,
\[
\frac{s_{n+1,u}}{t_{n+1,u}}=\left[a_{1},a_{2},\ldots,a_{m-1},a_{m}+1\right]
\quad \textrm{and} \quad \frac{s_{n+1,v}}{t_{n+1,v}}=\left[a_{1},a_{2},
\ldots,a_{m-1},a_{m}-1,2\right].\]
\end{lem}
\begin{proof}
For the first part of the lemma note that the second equality follows
by definition of $\mathfrak{T}_{n}$. The first equality is obtained
by induction as follows. We clearly have $\left\{
\left[2\right]\right\} =\mathfrak{T}_{1}\setminus\mathfrak{T}_{0}$.
Then assume that the assertion holds for $n-1$. Since the sets
$\mathcal{\mathfrak{T}}_{n}$
are $S$--invariant it follows for $n\geq 3$,\[
\mathfrak{T}_{n-1}\setminus\mathfrak{T}_{n-2}=\bigcup_{x\in\mathfrak{T}_{n-2}\setminus\mathfrak{T}_{n-3}}A(x)\cup
BS(x).\]
    For $\left[a_{1},\ldots,a_{k}\right]
\in\mathfrak{T}_{n-2}\setminus\mathfrak{T}_{n-3}$ we have by the
inductive assumption that $
\sum_{i=1}^k a_i =n-1$, and hence
   \begin{eqnarray*}
A\left(\left[a_{1},\ldots,a_{k}\right]\right) & = &
\frac{1}{1/\left[a_{1},\ldots,a_{k}\right]+1}=\left[a_{1}+1,a_{2},\ldots,a_{k}\right]\in
A_{k}^{n},\\
BS\left(\left[a_{1},\ldots,a_{k}\right]\right) & = &
\frac{1}{1+\left[a_{1},\ldots,a_{k}\right]}=\left[1,a_{1},a_{2},\ldots,a_{k}\right]\in
A_{k+1}^{n}.\end{eqnarray*}
   By combining the two latter observation, we obtain
\[\mathfrak{T}_{n-1}\setminus\mathfrak{T}_{n-2}\subset\bigcup_{k=1}^{n-1}
\bigcup_{(a_{1},a_{2},\ldots,a_{k}) \in A_{k}^{n}}\left[a_{1},a_{2},\ldots,a_{k}\right].\]
Since \begin{eqnarray*}
\card\left(\mathcal{\mathfrak{T}}_{n-1}\setminus\mathfrak{T}_{n-2}\right)
& = &
\card\left(\mathfrak{T}_{n-1}\right)-\card\left(\mathfrak{T}_{n-2}\right)=2^{n-2}\\
   & = & \sum_{k=1}^{n-1}{\binom{n-2}{
k-1}}=\card\left(\bigcup_{k=1}^{n-1}A_{k}^{n}\right),\end{eqnarray*}
   the first part of the lemma follows.

For the second part note that by the above
\[\left[a_{1},a_{2},\ldots,a_{m}+1\right],\,
\left[a_{1},a_{2},\ldots,a_{m}-1,2\right]\in
\mathfrak{T}_{n+1}\setminus\mathfrak{T}_{n}.\]
Since $\left[a_{1},a_{2},\ldots,a_{m}+1\right]$,
$\left[a_{1},a_{2},\ldots,a_{m}\right]$, and
$\left[a_{1},a_{2},\ldots,a_{m}-1,2\right]$ are consecutive neighbours
in $\mathfrak{T}_{n+1}$, the lemma follows.
\end{proof}

\begin{rem}
Note that $P$ can be written alternatively also in terms of denominators
of approximants as follows.
\[
P(\theta)=\lim_{n\to\infty}\frac{1}{n}\log\sum_{k=1}^{n}\sum_{\left(a_{1},\ldots,a_{k}\right)\in
A_{k}^{n}}q_{k}\left(\left[a_{1},\ldots,a_{k}\right]\right)^{-2\theta}.\]
In order to see this, note that  for $\theta\leq0$,
\[
\sum_{k=1}^{2^{n}}\left(t_{n,k}t_{n,k+1}\right)^{-\theta}\leq2\sum_{k=1}^{2^{n-1}}
\left(t_{n,2k}\right)^{-2\theta}\leq\sum_{k=1}^{2^{n+1}}\left(t_{n+1,k}t_{n+1,k+1}\right)^{-\theta}.\]
On the other hand, using the recursive definition
of $t_{n,k}$, we have for $\theta>0$,\[
\sum_{k=1}^{2^{n-1}}\left(t_{n-1,k}t_{n-1,k+1}\right)^{-\theta}\geq\sum_{k=1}^{2^{n-1}}
\left(t_{n,2k}\right)^{-2\theta}
\geq\frac{\left(n+1\right)^{-\theta}}{4}\sum_{k=1}^{2^{n+1}}\left(t_{n+1,k}t_{n+1,k+1}\right)^{-\theta}.\]
Therefore, by taking logarithms, dividing by $n$ and letting $n$
tend to infinity, we obtain  \[
P\left(\theta\right)=\lim_{n\to\infty}\frac{1}{n}\log\sum_{k=1}^{2^{n-1}}\left(t_{n,2k}\right)^{-2\theta}.\]
Hence, using Lemma \ref{lem:TnbyCF}, the result follows.
\end{rem}

\section{Dynamical cocycles versus Stern--Brocot
sequences\label{sec:Dynamical-cocycles-versus}}
In this section we introduce the dynamical cocycles which will be
crucial in the multifractal analysis to come. We show that 
 these cocycles are closely related to Stern--Brocot intervals and 
  continued fractions.  Finally,  we give the proof 
of the first part of Proposition \ref{Thm:2}. We remark that the results 
in this section could be obtained alternatively  by  using elementary 
estimates for countinued fractions only. Instead,  we have put  some
emphasis on obtaining these results by making use of the hyperbolic metric $d$ 
on ${\mathbb{H}}$. The intension here is that this should 
make it easier to follow the later transfer of the results 
of \cite{KesseboehmerStratmann:04a}, which were derived in terms of Kleinian 
groups, into the language 
of Stern--Brocot intervals and continued fractions.

Recall that the Poisson kernel $\PP$ for the upper half-plane is
given by\[
\PP:(z,\xi)\mapsto\frac{\Im\left(z\right)}{\left(\Re\left(z\right)-
\xi\right)^{2}+\Im\left(z\right)^{2}}  \,  ,  \textrm{  for all  } 
z\in\H, \xi\in\R.\]
   With $z_{0}$  defined as in Section \ref{sec:coding}, the cocycle
$I:\Sigma\to \R$ associated with the finite alphabet is given
by\[
I(x):=\left|\log\left(\PP\left(x_{1}(z_{0}),\pi\left(x\right)\right)\right)
-\log\left(\PP\left(z_{0},\pi\left(x\right)\right)\right)\right|   ,
 \textrm{  for  }  x=\left(x_{1},x_{2},\ldots\right)\in\Sigma.\]
We remark that $I$ is continuous with respect to the standard metric. 
Also, it is well--known that
$S_{n}I(x):=\sum_{i=0}^{n-1}I\left(\sigma^{i}\left(x\right)\right)$
is equal to the hyperbolic distance of $z_{0}$ to the horocycle
through $x_{1}x_{2}\cdots x_{n}\left(z_{0}\right)$ based at $\pi(x)$.
Furthermore, note that  in terms of the theory of iterations of maps,
$I$ is equal to 
the logarithm of the modulus of the derivative of the `twisted Farey 
map',  mentioned in  Remark \ref{rem:FG}. 

Similar,
we define the cocycle $I^{*}:\Sigma^{*}\to \R$ associated with the infinite
alphabet  as follows. For
$y=\left(X^{n_{1}},Y^{n_{2}},\ldots\right)\in\Sigma^{*}$ such that
$\left\{ X,Y\right\} =\left\{ A,B\right\} $, let $I^{*}$ be given by\[
I^{*}(y):=\left|\log\left(\PP\left(X^{n_{1}}Y(z_{0}),\pi^{*}\left(y\right)\right)\right)
-\log\left(\PP\left(z_{0},\pi^{*}\left(y\right)\right)\right)\right|.\]
One immediately verifies that 
$S_{k}I^{*}(y):=\sum_{i=0}^{k-1}I^{*}\left(\left(\sigma^{*}\right)^{i}(y)\right)$
is equal to the  hyperbolic distance of $z_{0}$ to the horocycle
based at $\pi^{*}(y)$ containing either the point $X^{n_{1}}Y^{n_{2}}\cdots
X^{n_{k}}Y\left(z_{0}\right)$
(if $k$ is odd) or $X^{n_{1}}Y^{n_{2}}\cdots Y^{n_{k}}X\left(z_{0}\right)$
(if $k$ is even). Note that  in terms of the theory of iterations of maps,
the function $I^{*}$ is clearly  an analogue of 
the logarithm of the modulus of the derivative of the Gauss map.
Throughout, we also require the potential function  $N:\Sigma^{*}\to\N$,
which is given by
\[ N\left((X^{n_{1}},Y^{n_{2}},\ldots)\right):=n_{1}, \hbox{ for each }
(X^{n_{1}},Y^{n_{2}},\ldots)\in \Sigma^{*}.\]
Note
that $S_{k}N((X^{n_{1}},Y^{n_{2}},\ldots))=\sum_{i=1}^{k}n_{i}$.

Finally,  the relevant potentials for the shift space $\left(\overline{\Sigma},\overline{\sigma}\right)$
are the functions 
\begin{equation} \label{barIN} \overline{I}:= I^{*}\circ
p_{A}=I^{*}\circ p_{B}\;\;\;\mbox{and}\;\;\;
\overline{N}:=N\circ p_{A}=N\circ p_{B}.\end{equation}
In here, 
$p_{X}$ refers to the inverse
branch with respect to $X\in\left\{ A,B\right\}$ of the $2$-$1$ factor map $p$ introduced in 
Section \ref{2.2}. More precisely, we have for 
$X,Y \in \left\{ A,B\right\} $ such that $X \neq Y$, \[
p_{X}\left(\left(n_{1},n_{2},n_{3},
\ldots\right)\right):=\left(X^{n_{1}},Y^{n_{2}},X^{n_{3}},\ldots\right).\]

The following lemma relates the Euclidean sizes of the Stern--Brocot 
intervals to the hyperbolic distances of $z_{0}$ to the elements in the orbit 
$G(z_{0})$. 

\begin{lem}
\label{lem:UniformlyExcursion} For each $n\in\N$ and $x\in\I$
such that $\pi^{-1}(x)=(x_{1},x_{2},...)\in\Sigma$, we have
   \[
\left|T_{n}(x)\right|\asymp m_{n}(x)\,
e^{-d\left(z_{0},x_{1}...x_{n}(z_{0})
\right)}.\]
In here, $m_{n}(x)$ is defined by $m_{n}(x):=\max\{ k:x_{n+1-i}=x_{n}\,\,\,\textrm{for}\,\,\,
i=1,...,k\}$.
\end{lem}
\begin{proof}
For $n=1$ the statement is trivial. For $n\geq2$, we first consider
the case $m_{n}(x)=1$. If $g:=x_{1}...x_{n-1}\in G$, then
   $g^{-1}(T_{n}(x))$ is equal to either $T_{1,1}$ (for
$x_{n}=A$) or $T_{1,2}$ (for $x_{n}=B$). Also, note that
for the modulus of the conformal derivative
   we have  \[
\left|\left(g^{-1}\right)'(\xi)\right|\asymp e^{d(z_{0},g(z_{0}))} ,
 \textrm{  for  } 
\xi \in T_{n}(x).\]
   Combining these two observations, we obtain \[
\left|T_{n}(x)\right|\asymp\left|g'|_{[0,1]}\right|\asymp\left|\left(g^{-1}\right)'|_{T_{n}(x)}\right|^{-1}\asymp
e^{-d(z_{0},g(z_{0}))}\asymp e^{-d(z_{0},gx_{n}(z_{0}))}.\]
   This proves the assertion for $m_{n}(x)=1$.

For the general situation we only consider the case
 $x_{1}\cdots x_{n}=A^{y_{1}}B^{y_{2}}\cdots B^{y_{k}}$. The
remaining cases can be dealt with in a similar way. In this case
$m_{n}(x)=y_{k}$,
and the above implies, for $l:=\sum_{i=1}^{k-1}y_{i}$,
\[
\left|T_{l+1}(x)\right|\asymp e^{-d\left(z_{0},x_{1}\cdots
x_{l+1}(z_{0})\right)}.\]
   Also, by using the well-known elementary fact that 
   $e^{d(z_{0},X^{k}(z_{0}))} \asymp \Im (X^{k}(z_{0})) \asymp 
   1/k^{2}$ for $X \in \{A,B\}$ and $k \in \N$, 
one immediately obtains  for $1<m \leq y_{k}$, \[
e^{d(z_{0},x_{1}\cdots x_{l+m}(z_{0}))}\asymp e^{d(z_{0},x_{1}\cdots
x_{l}(z_{0}))}e^{d(x_{1}\cdots x_{l}(z_{0}),x_{1}\cdots
x_{l+m}(z_{0}))}\asymp m^{2}e^{d(z_{0},x_{1}\cdots x_{l+1}(z_{0}))}.\]
Finally, one also immediately verifies that for $1\leq m \leq y_{k}$,
\begin{equation}
\left|T_{l+m}(x)\right|\asymp\sum_{k=m}^{\infty}k^{-2}\left|T_{l+1}(x)\right|\asymp
m^{-1}\left|T_{l+1}(x)\right|.\label{eq:AsympParabExcursion}\end{equation}
Combining the three latter observations, the statement of the lemma
follows.
\end{proof}

The previous lemma has the following immediate implication.
\begin{cor}
\label{cor:VergleichTSnI} For each $n\in\N$ and $x\in\I\,$
such that $\pi^{-1}(x) \in\Sigma$, we have
\[
\left|S_{n}I(\pi^{-1}(x))+\log\left|T_{n}(x)\right|\right|\ll\log n.\]
\end{cor}

The following lemma relates the cocycle $I^{*}$ to the 
sizes of the Stern--Brocot intervals and to the denominators $q_{k}^{2}$ of
the approximants.

\begin{lem}
\label{lem:Vergleichbar} For each $k\in\N$ and $x\in\I$
   we have, with $n_{k}(x):=S_{k}N\left(\left(\pi^{*}\right)^{-1}(x)
\right)$,
   \[
\left|T_{n_{k}(x) +1}(x)\right|\asymp\exp\left(-S_{k}I^{*}\left(\left(\pi^{*}\right)^{-1}(x)\right)\right)\asymp
q_{k}(x)^{-2}.\]
\end{lem}
\begin{proof}
We only consider the case $k$ even and $X=A$. The remaining
cases can be obtained in a similar way. Let $g:=A^{y_{1}}B^{y_{2}}\cdots
A^{y_{k}}\in G$, and note that then
\[
q_{k}(x)^{-2}\asymp e^{-d(z_{0},g(z_{0}))}.\]
Combining this with the fact that for $\xi\in T_{n+1}(x)$ we have
\[
\exp(-d(z_{0},g(z_{0})))\asymp\exp\left(-S_{k}I^{*}\left(\left(\pi^{*}\right)^{-1}(\xi)\right)\right)\]
(which follows since on $T_{n+1}(x)$
we have that
$\exp\left(S_{k}I^{*}\circ\left(\pi^{*}\right)^{-1}\right)$
is comparable to $|\left(g^{-1}\right)'|$), we obtain \[
e^{-S_{k}I^{*}\left(\left(\pi^{*}\right)^{-1}(x)\right)}\asymp
q_{k}(x)^{-2}.\]
   Finally, note that by Lemma \ref{lem:UniformlyExcursion} and since
$\exp\left(d\left(z_{0},gB(z_{0})\right)\right)\asymp\exp\left(d\left(z_{0},g(z_{0})\right)\right)$,
we have \[
\left|T_{n+1}(x)\right|\asymp e^{-d(z_{0},gB(z_{0}))}\asymp
e^{-d(z_{0},g(z_{0}))}.\]
   Combining these estimates, the lemma follows.
\end{proof}

We are now in the position to prove the first part of Proposition
\ref{Thm:2}.

\begin{proof}
[Proof of first part of Proposition  \ref{Thm:2}] The equalities
$\ell_{3}=\ell_{6}$ and $\ell_{1}=\ell_{5}$
 are immediate consequences of the following well--known Diophantine
inequalities. For all
$x\in[0,1]$ and $k\in\N$, we have (see e.g. \cite{Khinchin:35})
  \begin{equation}
\frac{1}{q_{k}(x)\left(q_{k+1}(x)+q_{k}(x)\right)}<\left|x-\frac{p_{k}(x)}{q_{k}\left(x\right)}
\right|<\frac{1}{q_{k}(x)q_{k+1}(x)}.\label{eq:CFInequality}\end{equation}

In order to show that $\ell_{1}=\ell_{4}$,  let $n_{k}(x):=S_{k}N
\left(\left(\pi^{*}\right)^{-1}(x)\right)$, for $x\in\I$ and $k\in\N$. 
Obviously, $\left(\log q_{k}(x)/n_{k}(x)\right)$ is a subsequence 
of the sequence $\left(-\log |T_{n}(x)| / n\right)$. Therefore, if 
$\ell_{4}(x)$ exists then so does $\ell_{1}(x)$, and both limits 
must coincide. For the reverse, suppose that $\ell_{1}(x)$ exists 
such that $\ell_{1}(x) = \alpha$. Let $m_{n}(x)$ be defined as in 
the statement of Lemma \ref{lem:UniformlyExcursion}, and put
$k_n(x):=\sup\{k\in \N:n_k(x)\leq n\}$. By combining (\ref{eq:AsympParabExcursion}) and Lemma
\ref{lem:Vergleichbar}, we then have
\begin{eqnarray*}
\lim_{n\to\infty}\frac{-\log\left|T_{n}\left(x\right)\right|}{n} & \leq &
 \lim_{n\to\infty}\frac{2\log q_{k_{n}\left(x\right)}\left(x\right)+
2\log (m_{n}\left(x\right)+1)}{n_{k_{n}\left(x\right)}\left(x\right)+m_{n}\left(x\right)}\\
 & \leq & \lim_{n\to\infty}\frac{2\log q_{k_{n}\left(x\right)}\left(x\right)}{n_{k_{n}\left(x\right)}\left(x\right)}
+\lim_{n\to\infty}\frac{2\log (m_{n}\left(x\right)+1)}{n_{k_{n}\left(x\right)}\left(x\right)+m_{n}\left(x\right)}=
\alpha.\end{eqnarray*}
This gives the upper bound, and hence finishes the proof 
in particular for $\alpha=0$. 
For the opposite inequality we can therefore assume without loss of 
generality that $\alpha>0$. First, observe that
\begin{eqnarray*}
\lim_{n\to\infty}\frac{-\log\left|T_{n}\left(x\right)\right|}{n} & = &
 \lim_{n\to\infty}\frac{2\log  q_{k_{n}\left(x\right)}\left(x\right)+
2\log m_{n}\left(x\right)}{n_{k_{n}\left(x\right)}\left(x\right)+m_{n}\left(x\right)}\\
 & \geq & \lim_{n\to\infty}\frac{2\log q_{k_{n}\left(x\right)+1} 
 \left(x\right)}{n_{k_{n}\left(x\right)+1}\left(x\right)}-
\lim_{n\to\infty}\frac{2\log 
a_{k_{n}\left(x\right)+1}\left(x\right)}{n_{k_{n}\left(x\right)}\left(x\right)} .
\end{eqnarray*}
Hence, it is now sufficient to show that $\lim_{k\to\infty}\log 
a_{k+1}\left(x\right) / n_{k}\left(x\right)=0$, 
or what is equivalent   $\lim_{k\to\infty} \log 
a_{k+1}\left(x\right) / \log q_{k}\left(x\right) =0$.
 For this, observe that
\[ \lim_{k\to \infty} \frac{\log q_{k+1}(x)}{n_{k+1}(x)} = 
\lim_{k\to \infty} \frac{\log q_{k}(x) +\log a_{k+1}(x)}{n_{k} (x)+ 
a_{k+1}(x)} =  
\lim_{k\to \infty} \frac{\log q_{k}(x)  \left(1+ \frac{\log 
a_{k+1}(x)}{\log q_{k}(x)} \right)}{n_{k}(x) \left(1+ 
\frac{a_{k+1}(x)}{n_{k}(x) } \right) } .\]
If we would have that  $\lim_{k\to\infty} \log 
a_{k+1}\left(x\right) / \log q_{k}\left(x\right) \neq 0$, then 
there exists a subsequence $(k_l)$ such that  $\lim_{l\to\infty} \log 
a_{k_l+1}\left(x\right) / \log q_{k_l}\left(x\right)= 
c$, for some $c \in (0,\infty]$.
It follows that $\lim_{l \to\infty} a_{k_l+1} = \infty $, and hence by 
combining this with  the 
calculation above, we obtain
\[
1=\lim_{l\to\infty}\frac{\log a_{k_{l}+1}(x)\cdot n_{k_l}(x)}{a_{k_{l}+1}(x) \cdot\log q_{k_l}(x)}=\frac{0}{\alpha}=0.
\]
This shows that $\ell_{1}(x)=\ell_{4}(x)$, and hence finishes the proof. 
\end{proof}

\section{\label{sec:The-shapes-of-P} Analytic properties of $P$ and
$\widehat{P}$}
The main goal in this section is to derive various analytic properties 
of  the
Stern--Brocot pressure function $P$. These properties are derived by 
considering the pressure functions associated with the systems $\Sigma, \Sigma^{*}$ 
and $\overline{\Sigma}$. In order to introduce these functions, let $\mathcal{C}_{n}:=\left\{ C_{n}(x):x\in\Sigma\right\} $ refer
to the set of all $n$--cylinders  \[
C_{n}\left(x\right):=\left\{ y\in\Sigma:y_{i}=x_{i}, \hbox{for } i=1,\ldots,n\right\}
.\]
Likewise, let $\mathcal{C}^{*}_{n}$ (resp. 
$\overline{\mathcal{C}}_{n}$) refer to the set of $n$-cylinders for 
the system $(\Sigma^{*},\sigma^{*})$ (resp. 
$(\overline{\Sigma},\overline{\sigma})$). The pressure function
$\mathcal{P}$ associated with $\Sigma$ is then  given by \[
\mathcal{P}(\theta):=\lim_{n\to\infty}\frac{1}{n}\log\sum_{C\in\mathcal{C}_{n}}
\exp\left(\sup_{x\in
C}S_{n}\left(-\theta I\right)(x)\right),  \textrm{ for }
\theta\in\R.\]
Also,  for the system $\Sigma^{*}$ we define the pressure functions 
$\P^{*}$ and $P^{*}$,  for $\theta<1, q>0$ and $f: \Sigma^{*} 
\to\R$ continuous,   by
\[
\P^{*}(f):=\lim_{n\to\infty}\frac{1}{n}\log\sum_{C^{*}\in\mathcal{C}_{n}^{*}}\exp\left(\sup_{y\in
C^{*}}S_{n}f\left(y\right)\right) \hbox{  and  } P^{*} (\theta, q):= \P^{*}(-\theta
I^{*}-qN).\]
Finally,  the pressure functions $\overline{\P}$ and $\overline{P}$
associated with $(\overline{\Sigma}, \overline{\sigma})$
are given completely analogous,  for $\theta<1, q>0$ and 
$g: \overline{\Sigma} \to\R$ continuous, by
\[ 
\overline{\P}(g):=\lim_{n\to\infty}\frac{1}{n}\log\sum_{\overline{C}\in\overline{\mathcal{C}}_{n}}\exp\left(\sup_{y\in
\overline{C}}S_{n}g\left(y\right)\right) 
\hbox{  and  }  \overline{P}(\theta, q):= \overline{\P}(-\theta
\overline{I}-q\overline{N}).\]
Clearly,  by recalling the definitions of $\overline I, \overline N$ 
and $I^{*}, N$ in Section 3,  we immediately have that $\overline{P} = P^*$.
\subsection{
Analytic Properties of 
$P^{*}$ by Hanus, Mauldin and Urba\'nski}
In this subsection we employ  important results of  Hanus, Mauldin and 
Urba\'nski obtained in \cite{MauldinUrban:01}  
and \cite{HanusMauldiUrban:02}.  The results here will be
crucial cornerstones  in our 
subsequent analysis of 
the Stern--Brocot pressure. 

Studies of analytic properties of 
 pressure functions are usually based on the existence
of certain Gibbs measures, here on $\Sigma^{*}$ and 
$\overline{\Sigma}$. The existence of these measures in our situation here
is guaranteed by the following proposition, which essentially follows 
from a result in   
\cite{MauldinUrban:01}.
\begin{prop}\label{propHMU2}
For each  $\theta <1$, $ q>0$, and for $(\theta,q)=(1,0)$, 
there exists a unique completely ergodic $\overline{\sigma}$--invariant Gibbs measure
$\overline{\mu}_{\theta,q}$ 
associated with the potential $-\theta \overline{I}-q\overline{N}$. 
That is,  we have  for all  $n\in\N, \overline{C}\in
\overline{\mathcal{C}}_{n}$ and $y\in \overline{C}$, 
\begin{equation}
	\overline{\mu}_{\theta,q}(\overline{C})\asymp\exp\left(S_{n}\left(-\theta
	\overline{I}(y)-q\overline{N}(y)\right)-n\,\overline{\P}(-\theta
	\overline{I}-q\overline{N})\right).\label{eq:Gibbs*Property}
\end{equation}
In particular,  the Borel measure
$\mu_{\theta,q}^{*}:=1/2\cdot \left(\overline{\mu}_{\theta,q}\circ
p_{A}^{-1}+\overline{\mu}_{\theta,q}\circ p_{B}^{-1}\right)$
is an ergodic  $\sigma^{*}$--invariant  Gibbs 
measure  on $\Sigma^{*}$ such that
for all $n \in \N, C^{*}\in \mathcal{C}_{n}^{*}$  and $y\in C^{*}$, 
\begin{equation}
\mu^{*}_{\theta,q}(C^{*})\asymp\exp\left(S_{n}\left(-\theta
I^{*}(y)-q N (y)\right)-n\, \P^{*}(-\theta
I^{*}-q N)\right).\label{eq:Gibbs*Property2}\end{equation}
The measure $\mu_{\theta,q}^{*}$ is unique with respect
to this property, and  $\overline{\mu}_{\theta,q}=\mu_{\theta,q}^{*}\circ
p^{-1}$.
\end{prop}
\begin{proof} By \cite{KesseboehmerStratmann:04a} (Lemma
3.4), the cocycle $I^{*}$ is H{\"o}lder continuous in the sense that there
exists $\kappa>0$ such that for each $n\in\N$, \[
\sup_{C\in\overline{\mathcal{C}}_{n}}\sup_{x,y\in
C}\left|\overline{I}(x)-\overline{I}(y)\right|\ll\exp(-\kappa\, n).\]
Clearly, we also immediately have that $\overline{N}$ is  H\"older 
continuous.
Furthermore,  the following summability condition holds for $\theta <1$, $ q>0$, and for $(\theta,q)=(1,0)$, 
\begin{equation}\label{summability}
\sum _{i\in \N}\exp(\sup\{-\theta \overline{I}(x)-q 
\overline{N}(x):x_1=i\})\ll\sum _{i\in \N}
i^{-2\theta}\cdot \e^{-qi}<\infty .
\end{equation}
Hence, all preliminaries of  \cite{MauldinUrban:01} (Corollary 2.10) 
are fulfilled, which then gives the existence of a 
unique invariant Gibbs measure $\overline{\mu}_{\theta,q}$ with properties as stated in the 
proposition.

Immediate consequences of the definition of $N$ and the definition $\mu_{\theta,q}^{*}:=1/2\cdot \left(\overline{\mu}_{\theta,q}\circ
p_{A}^{-1}+\overline{\mu}_{\theta,q}\circ p_{B}^{-1}\right)$ are that  
$\mu_{\theta,q}^{*}$
is $\sigma^*$--invariant, that $\mu_{\theta,q}^{*}$ fulfills the Gibbs property (\ref{eq:Gibbs*Property2}), and 
that the equality 
$\overline{\mu}_{\theta,q}=\mu_{\theta,q}^{*}\circ
p^{-1}$ is satisfied. To prove  
ergodicity of  $\mu^*_{\theta,q}$, 
let $D\subset\Sigma^*$ such that ${\sigma^*}^{-1}(D)=D$.
We then have $\overline{\sigma}^{-1} p_X^{-1}(D)=p_Y ^{-1}(D)$,
for $X,Y\in \{A,B\}$ such that $X\neq Y$. This 
gives $\overline{\sigma}^{-2}({p_X}^{-1}(D))={p_X}^{-1}(D)$.
Since $\overline{\mu}_{\theta,q}$ is completely ergodic, which by 
definition means that  $\overline{\mu}_{\theta,q}$ is ergodic with 
respect to $\overline{\sigma}^n$ for all $n\in \N$, it follows 
 $\overline{\mu}_{\theta,q}({p_X}^{-1}(D))\in\{0,1\}$. The $\overline{\sigma}$-invariance of 
 $\overline{\mu}_{\theta,q}$ then implies that $\overline{\mu}_{\theta,q}({p_X}^{-1}(D))=
\overline{\mu}_{\theta,q}(\overline{\sigma}^{-1}({p_X}^{-1}(D)))=\overline{\mu}_{\theta,q}({p_Y}^{-1}(D))$.
Consequently, it follows that $\mu^*_{\theta,q}(D)\in \{0,1\}$.
\end{proof}

 The following proposition employs yet  another result of Hanus, Mauldin and 
 Urba\'nski, obtained in their spectral analysis of the  Perron--Frobenius operator.
\begin{prop}\label{propHMU3}
    The pressure function $P^{*}$ is a convex, decreasing and real-analytic 
    function with
    respect to both coordinates  In the second coordinate $P^{*}$ is strictly
    decreasing to $(-\infty)$.  In particular, there hence exists a positive
    real--analytic function $\beta$ on $(-\infty,1)$ such that
    $P^{*}(\theta,\beta(\theta))=0$.  Furthermore, for the
    derivative of $\beta$ at $\theta<1$ we have 
    \begin{equation}\label{derive} \beta'(\theta)=\frac{-\int I^{*}\,
    d\mu_{\theta}^{*}}{\int N\, d\mu_{\theta}^{*}}=-\int I\, d\mu_{\theta}.
    \end{equation}
   In here,
    $\mu_{\theta}^{*}:=\mu_{\theta,\beta(\theta)}^{*}$ refers to the
    unique $\sigma^{*}$--invariant Gibbs measure  
    associated with the
    potential $-\theta I^{*}-\beta(\theta)N$. Also, $\mu_{\theta}$ 
    refers to the $\sigma$--invariant probability measure on $\Sigma$
    absolutely continuous to  $\mu_{\theta}^{*}$, whose existence is 
    guaranteed by {\em Kac's formula} (\cite{Kac:47}; see the proof).
\end{prop}
\begin{proof} First, note that it is sufficient to verify the 
statements in the proposition for $\overline{\Sigma}$ only. Then note 
that
 $\overline{I}$ and $\overline{N}$  are H\"older continuous, that
the summability condition (\ref{summability}) is  satisfied
for all 
$(\theta,q)\in (-\infty,1)\times (0,\infty)$, and that
\[
\int \left( \theta  \overline{I} + q \overline{N} \right) 
d\overline{\mu}_{\theta,q}\ll\sum_{n\in\N}(2\theta\log n+q n) n^{-2\theta}\e ^{-q n}<\infty.
\]
 Hence, we can apply \cite{HanusMauldiUrban:02}  (Proposition 6.5)  
(see also \cite{MauldinUrbanski:03} (Proposition 2.6.13)), from which 
it follows that $\overline{P}(\theta,q)$ is real-analytic on 
$(-\infty,1)\times (0,\infty)$. Next, note that for the partial derivatives of 
$\overline{P}$ 
we have
\[
\frac{\partial \overline{P}(\theta,q)}{\partial \theta}=\int -\overline{I} \;d 
\overline{\mu}_{\theta,q} \hbox{  and  }  
\frac{\partial \overline{P}(\theta,q)}{\partial q}=\int -\overline{N} \;d \overline{\mu}_{\theta,q}<0.
\]  
This shows that  $\overline{P}$ as a function  in the second 
coordinate is strictly decreasing,  which then gives the existence of 
a 
real-analytic function $\beta:(-\infty,1)\to (0,\infty)$ for which 
$\overline{P}(\theta,\beta(\theta))=0$, for all $\theta <1$. 
Now, the first equality in (\ref{derive}) follows
from the Implicit Function Theorem. The second equality is a consequence
of {\em Kac's formula}  (\cite{Kac:47}),  which guarantees that there exists a 
$\sigma$--invariant  measure $\widetilde{\mu}_{\theta}$
on $\Sigma$,
given by \begin{equation}
\widetilde{\mu}_{\theta}(M):=\int\sum_{i=0}^{N(y)-1}\1_{M}\circ\sigma^{i}(\iota(y))\,\,
d\mu_{\theta}^{*}(y), \hbox{ for each Borel set } M \subset \Sigma .\label{eq:Kac}
\end{equation} In here,
$\iota:\Sigma^{*}\to\Sigma$ refers to the canonical injection which
maps an element of $\Sigma^{*}$ to its representation by means of the
finite alphabet of $\Sigma$.  We remark that Kac's formula gives in 
fact a 1--1 correspondence between $\sigma$--invariant measures on 
$\Sigma$ and $\sigma^{*}$--invariant measures on $\Sigma^{*}$.
Now, one immediately verifies that \[
   \widetilde{\mu}_{\theta}(\Sigma)=\sum_{\ell=1}^{\infty} \ell \mu^{*}_{\theta}\left(\left\{
   N=\ell\right\}
   \right)\asymp\sum_{\ell=1}^{\infty}\ell^{-2\theta+1}e^{-\beta(\theta)\ell}<\infty.\]
   Hence, this allows to define 
   $\mu_{\theta}:=\widetilde{\mu}_{\theta}/\widetilde{\mu}_{\theta}
   \left(\Sigma\right)$.
   In particular, we then have 
   that $\int I\,
   d\mu_{\theta}=\mu_{\theta}^{*}(N)^{-1}\int I^{*}\, 
   d\mu_{\theta}^{*}$, from which the second equality follows.
Finally, note that since $\overline{P}$
 and $P^*$ coincide,  in the integrals above we can replace 
 $\overline{I}$ and $\overline{N}$ by $I^{*}$ and $N$. This finishes proof 
 of the proposition.
 \end{proof}
 
 \begin{rem}
     Note that the measure $\mu_{\theta}$ in Proposition \ref{propHMU3}  
     is in fact a weak Gibbs measure for 
       the potential $-\theta I$. Therefore,  the results of
       \cite{Kesseboehmer:01} are applicable, and hence in this way 
       one could immediately obtain some Large Deviation results for the 
       situation here.
     \end{rem}

\subsection{Analytic properties of $P$ and
$\widehat{P}$}
In this subsection we employ the results of  the previous subsection, 
in order to derive analytic properties 
of the Stern--Brocot pressure
$P$ and its Legendre transform $\widehat{P}$. 

A key 
preliminary observation is stated in the following proposition, which 
shows that the Stern--Brocot pressure $P$ coincides with the function 
$\beta$ obtained in 
Proposition \ref{propHMU3}.

\begin{prop}\label{prop:P=beta}
 For $\theta <1$, we have
 \[ P(\theta)=\beta(\theta).\]
 \end{prop}
 \begin{proof}
   Let $\mu_{\theta}^{*}$ be the measure obtained in Proposition 
   \ref{propHMU3}, for $\theta<1$ fixed. 
   First, recall from the proof of Proposition  \ref{propHMU3} 
   that the measure class of  $\mu^{*}_{\theta}$ contains a 
$\sigma$--invariant  probability measure $\mu_{\theta}$
on $\Sigma$ for which  $\int I\,
   d\mu_{\theta}=\mu_{\theta}^{*}(N)^{-1}\int I^{*}\, 
   d\mu_{\theta}^{*}$. Secondly, note that for the measure $\overline{\mu}_{\theta}:=
   \mu_{\theta}^{*}\circ
   p^{-1}$ we have by  
   {\em Abramov's formula} (\cite{Abramov:59}, \cite{Neveu:69}) that
   the 
   measure theoretical entropies  $h_{\mu_{\theta}}$ and 
   $h_{\overline{\mu}_{\theta}}$ are related by $h_{\mu_{\theta}} = h_{\overline{\mu}_{\theta}}
   /\mu_{\theta}^{*}(N)$. Thirdly,  
   by applying {\em  Pinsker's result on relative 
   entropies} (\cite{Rohlin:67})  to our situation here, we obtain  that the relative entropy
$h_{\mu^{*}_{\theta}}\left(\sigma^{*}|\overline{\sigma}\right)$
of $\mu^{*}_{\theta}$ vanishes. This gives
$h_{\mu^{*}_{\theta}}-h_{\overline{\mu}_{\theta}}=
h_{\mu^{*}_{\theta}}\left(\sigma^{*}|\overline{\sigma}\right)=0$.
Combining these observations with the usual variational principle 
(\cite{DenkerGrillenbergerSigmund:76}), 
   it now follows
   \begin{eqnarray*}
   P(\theta) & \geq & h_{\mu_{\theta}}-\int  \theta I\,
   d\mu_{\theta}
       = 
   \left(\mu_{\theta}^{*}(N)\right)^{-1}\left(h_{\mu_{\theta}^{*}}-\int\theta
   I^{*}\, d\mu_{\theta}^{*}\right)\\
      & = &
   \left(\overline{\mu}_{\theta}(\overline{N})\right)^{-1}
   \left(h_{\overline{\mu}_{\theta}}-\int\theta
   \overline{I} \,
   d\overline{\mu}_{\theta}\right)
       =  \beta(\theta).\end{eqnarray*}
   In here, the latter equality is obtained as follows. Note  that
   $\overline{\mu}_{\theta}$
   is an equilibrium measure on $\left(\overline{\Sigma},\overline{\sigma}\right)$
   for the potential $-\theta \overline{I}-\beta(\theta)\overline{N}$.
   Also, $\overline{\mu}_{\theta}$ fulfills the Gibbs property  
   \ref{eq:Gibbs*Property}, and 
   $\overline{\mu}_{\theta}\left(\theta
   \overline{I}+\beta(\theta)\overline{N}\right) < \infty$. Next, recall
   \emph{Sarig's variational principle} (\cite{Sarig:99})  which 
   states that
   for 
   $g: \overline{\Sigma} \to\R$ H\"older continuous, 
   \begin{equation}
   \overline{\P}(g)=\sup\left\{ h_{\overline{\mu}} +\int g\,
   d\overline{\mu} :\overline{\mu}\in\mathcal{M}\left(\overline{\Sigma} 
   ,\overline{\sigma}\right)\,\,\textrm{such
   that}\,\, -\int g\, d\mu < \infty\right\} .\label{eq:Sarig}\end{equation}
 In here, 
$\mathcal{M}(\overline{\Sigma},\overline{\sigma})$ refers to the set of
$\overline{\sigma}$--invariant Borel probability measures on 
$\overline{\Sigma}$.  Applying this to the situation here, we obtain
  \[ h_{\overline{\mu}_{\theta}}-\int \left(\theta
   \overline{I} + \beta(\theta) \overline{N} \right) \,
   d\overline{\mu}_{\theta} = \overline{\P} (-\theta\overline{I} - 
   \beta(\theta) \overline{N})=  \overline{P}(\theta, \beta(\theta)) = 
   0.\]
   An elementary rearrangement  then gives the result.

   For the reverse inequality, 
   first note that 
   we can induce $\left(\Sigma,\sigma\right)$ on the set 
   $\mathcal{H}:=\left\{
      \left(x_{1},x_{2},\ldots\right)\in\Sigma:x_{1}\not=x_{2}\right\} $.
      For the resulting induced system  the return 
      time to $\mathcal{H}$ of a point
   $y=\iota(X,Y^{n_{1}},X^{n_{2}},\ldots)\in\mathcal{H}$ is given by
   $n_{1}=N\left(\sigma^{*}(X,Y^{n_{1}},X^{n_{2}},\ldots)\right)$. 
   Define 
   $\mathcal{G}:=\iota\left(\Sigma^{*}\right)\cap\mathcal{H}$, and let 
   $m_{\theta}\in\mathcal{M}(\Sigma,\sigma)$
   be an ergodic equilibrium measure for the potential $-\theta I$, that
   is $P(\theta)=h_{m_{\theta}}-\theta\int I\, dm_{\theta}$.
   In this situation we necessarily have that 
   $m_{\theta}(\mathcal{G})>0$, and  this can be seen as follows.
   First, we show that  $m_{\theta}(\mathcal{G})=0$ implies that
    $m_\theta$ is  equal to either $\delta_{\overline{A}}$ or 
    $\delta_{\overline{B}}$, where $\delta_{\overline{A}}$ (resp.  $\delta_{\overline{B}}$)  refers to
the Dirac measure at the periodic point $\overline{A}=\pi^{-1}(0)$ (resp. $\overline{B}=\pi^{-1}(1)$). 
   Namely, if $m_\theta(\mathcal{H})=0$ then 
    we immediately have $m_\theta\in\{\delta_{\overline{A}}, 
    \delta_{\overline{B}}\}$. On the other hand, if  $m_\theta(\mathcal{H})>0$ then
    ergodicity  of $m_{\theta}$ gives $m_\theta(\Sigma\setminus\iota(\Sigma^*))=1$. 
    Now, since $I(x)\geq 0$ for all $x \in \Sigma$,  where $I(x) = 0$ if 
    and only if  
    $x\in \{\overline{A},\overline{B}\}$, we have 
     $\lim n^{-1}S_n I(x)=0=\int I \;dm_\theta$  for $m_\theta$-almost 
     every $x$, and this again implies that $m_\theta\in\{\delta_{\overline{A}},
      \delta_{\overline{B}}\}$. Thus, if 
      $m_{\theta}(\mathcal{G})=0$ then $m_\theta\in\{\delta_{\overline{A}},
      \delta_{\overline{B}}\}$. This shows that
   $h_{m_{\theta}}= $ $m_{\theta}\left(-\theta I\right)=0$, giving 
   $P(\theta)=0$, and hence contradicting the fact 
   $P(\theta)\geq\beta(\theta)>0$. Therefore, we can assume without loss 
   of generality that $m_{\theta}(\mathcal{G})>0$. We can now use  Kac's formula
   once more,  
   which
     guarantees 
   that there exists  a $\sigma^{*}$--invariant probability measure $m_\theta^{*}$
   in the 
   measure class of $m_{\theta}$, such that  $m_\theta^{*}:=
   \frac{1}{m_\theta\left(\mathcal{G}\right)}
   m_\theta|_{\mathcal{G}}\circ\sigma^{-1}\circ\iota$  and
   \[
   -\int\left(-\theta I^{*}-\beta(\theta)N\right)\,
   dm_{\theta}^{*}=\left(m_{\theta}(\mathcal{G})\right)^{-1}\left(\int\theta
   I\, dm_{\theta}+\beta(\theta)\right)<\infty.\]
      For $\overline{m}_{\theta}:=m_{\theta}^{*}\circ p^{-1}$,  we  
      argue similar as above and obtain

   \begin{eqnarray*}
   0 & \geq & h_{\overline{m}_{\theta}}-\int\left(\theta
   \overline{I}+\beta(\theta)\overline{N}\right)\,
   d\overline{m}_{\theta} \,\,\, \qquad\quad(\emph{Sarig's 
   variational principle})\\
    & = & h_{m_{\theta}^{*}}-\int\left(\theta
   I^{*}+\beta(\theta)N\right)\,
   dm_{\theta}^{*}\,\quad\qquad(\emph{Pinsker's result on relative 
   entropies})\\
      & = & m_{\theta}^{*}(N)\left(h_{m_{\theta}}-\int\theta I\,
   dm_{\theta}-\beta(\theta)\right) \qquad(\emph{Abramov's
   formula})\\
      & = &
   m_{\theta}^{*}(N)\left(P(\theta)-\beta(\theta)\right)
   \,\,\,\,\,\,\,\qquad\quad\,\,\,\,\,\,\,\,\,\,\,\,\,\,
   \,(\emph{since} \,\, m_{\theta}\, \,\emph{is an equilibrium
   state}).\end{eqnarray*}    
   \end{proof}
The following proposition  collects  the properties
of $P$ and $\widehat{P}$ which will be crucial in the analysis to come.
\newpage
\begin{prop}\label{pro:AnalyticPropertiesP}~\nopagebreak
\begin{enumerate}
\item The Stern--Brocot pressure function $P$ coincides with the
pressure function
$\mathcal{P}$ associated with $\Sigma$.  
\item $P$ is convex and non-increasing on $\R$ and real-analytic on
$(-\infty,1)$.
\item $P(\theta)=0$, for all $\theta\geq1$.
\item $P$ is differentiable throughout $\R$.
\item The domain of $\widehat{P}$ is equal to $[-\alpha_{+},0]$, where
\[
-\alpha_{+}:=\lim_{\theta \to-\infty}\frac{P(\theta)}{\theta}=-2\log\gamma.\]
\item We have
$\lim_{\alpha\searrow0} \widehat{P}\left(-\alpha\right)/(-\alpha)=1$.
\item We have
$\lim_{\alpha\nearrow2\log\gamma}\left(-\widehat{P}\left(-\alpha\right)\right)=0.$
\item We have $\lim_{\theta\to-\infty}\left(P(\theta)
+2\theta\log\gamma\right)=0$.
\end{enumerate}
\end{prop}
For the proofs of (7) and (8) the following lemma will turn out
to be useful.
\begin{lem}
\label{lem:QGammaInequalitystrict} For each
$x:=[a_{1},a_{2},a_{3},\ldots]
\in(0,1)$
and $k \in\N_0$ we have, with $\tau_0:=0$, $\tau_{k}:=\sum_{i=1}^{k}a_{i}$ for 
$k\in \N$,
and $\rho:=1-\gamma^{-6}$, \[
q_{k}(x)\leq\gamma^{\tau_{k}}\rho^{\tau_{k}-k-1}.\]
\end{lem}


\begin{proof}
We give a  proof by complete induction of the slightly stronger inequality
\begin{equation}
q_{k}(x)\leq\gamma^{\tau_{k}}\rho^{\tau_{k}-k} \rho^{\delta_{1 ,
a_{k}}-1},\label{eq:slightlyStronger}\end{equation} in which
$\delta$ denotes the Kronecker symbol.  \\
First note that $q_0\equiv 1$,  $q_{1}([1,\ldots])=1\leq\gamma^{1}\rho^{1-1}$, and
if $a_{1}\geq2$ then one immediately verifies that
$q_{1}\left[a_{1},\ldots\right]=a_{1}\leq\gamma^{a_{1}}\rho^{a_{1}-1}\rho^{-1}$.
Also, for $k \in \N$ we have \begin{equation} q_k(\gamma-1)=q_k([1,1,1,\ldots])=f_k\leq \gamma^{k}=\gamma
^{\tau_k}\rho^{\tau_k-k},\label{eq:abschGoldenMean}\end{equation}
where $f_{k}$ 
 refers to the $(k+1)$-th member of the Fibonacci sequence 
 $(f_{0},f_{1},f_{2},\ldots):= (0,1,1,2,\ldots)$, given by 
 $f_{k+1}:= f_{k-1}+f_{k}$ for all $k \in \N$.
 Recall that  $f_{k}=\left(\gamma^{k}-\left(-\gamma\right)^{-k}\right)/\sqrt{5}$.
Now suppose that (\ref{eq:slightlyStronger}) holds for some $k\in \N$ and
  for all $0\leq m\leq k$. It is then sufficient to consider the
following
two cases.
\begin{enumerate}
\item \label{enu:case1} If $a_{k+1}=1$ such that $a_{n}\geq2$ and $a_{n+i}=1$, for all
$i=1,\ldots,l$
and some $n\leq k$ and $l\geq k-n+1$, then
$q_{n-1}(x)\leq\gamma^{\tau_{n-1}}\rho^{\tau_{n-1}-n+1} \rho^{-1}$
and
$q_{n}(x)\leq\gamma^{\tau_{n}}\rho^{\tau_{n}-n} \rho^{-1}$.  Hence,
 an elementary
calculation gives \begin{eqnarray*} q_{n+l}(x) & = &
f_{l+1}q_{n}(x)+f_{l}q_{n-1}(x)\\
   & \leq &
f_{l+1}\gamma^{\tau_{n}}\rho^{\tau_{n}-n} \rho^{-1}+f_{l}\gamma^{\tau_{n-1}}
\rho^{\tau_{n-1}-n+1} \rho^{-1} \\
   & \leq &
\gamma^{\tau_{n+l}}\rho^{\tau_{n+l}-n-l}\left(\rho^{-1}\left(\frac{f_{l+1}}
{\gamma^{l}}+\frac{f_{l}}
{\gamma^{a_{n}+l}\rho^{a_{n}-1}}\right)\right)\\
   & \leq &
\gamma^{\tau_{n+l}}\rho^{\tau_{n+l}-n-l}\underbrace{\left(\rho^{-1}
\left(\frac{f_{l+1}}{\gamma^{l}}+
\frac{f_{l}}{\gamma^{l}\left(\gamma\rho\right)^{2}}\right)\right)}_{\leq1}.
\end{eqnarray*}

\item If $a_{k+1}=2$, then either  $a_i=1$ for $i=1,\ldots,k$, or
there exists  $n\leq k$ such that  $a_{n}\geq2$ and $a_{i}=1$ for all $i$ with
$n <i\leq k$. In the first case we use (\ref{eq:abschGoldenMean}), whereas in 
the second case we employ (\ref{enu:case1}), and obtain
  \begin{eqnarray*}
q_{k+1}\left(\left[a_{1},\ldots,a_{k},2\right]\right) & = &
q_{k+2}\left(\left[a_{1},\ldots,a_{k},1,1\right]\right)\\
   & \leq &
\gamma^{\tau_{k+1}}\rho^{\tau_{k+1}-k-1} \rho^{-1}.\end{eqnarray*}
For $a_{k+1}>2$, the inequality follows by induction over $a_{k+1}$,
using  (\ref{enu:case1})  and the fact that
$q_{k+1}\left(\left[a_{1},\ldots,a_{k},a_{k+1}\right]\right)=q_{k+2}
\left(\left[a_{1},
\ldots,a_{k+1}-1,1\right]\right)$.

\end{enumerate}
\end{proof}
Before giving the proof of Proposition \ref{pro:AnalyticPropertiesP},
we remark that the statements  (7) and (8)  in Proposition \ref{pro:AnalyticPropertiesP}
are in fact equivalent.
Nevertheless, we shall prove these two statements separately, where
the proof of (7) primarily uses ergodic theory, whereas the
proof of (8) is of elementary number theoretical nature.

\begin{proof}
[Proof of Proposition \ref{pro:AnalyticPropertiesP}]~

\textbf{\emph{ad (1)}.} \, The assertion is an immediate consequence of 
(\ref{eq:AsympParabExcursion}) and
Corollary \ref{cor:VergleichTSnI}.

\textbf{\emph{ad (2)}.} \, The assertion follows immediately by combining 
Proposition \ref{propHMU3} and Proposition \ref{prop:P=beta}. 
Alternatively, the statement can also be derived from 
Proposition 2.1 in
\cite{KesseboehmerStratmann:04a}.

\textbf{\emph{ad (3)}.}  \, By definition of $P$ we have $P(1)=0$. 
Also, 
by (2)  we know that $P$ is non-increasing. Therefore, it is sufficient to show that $P$
is non-negative, and indeed this follows  since\[
P(\theta)=\lim_{n\to\infty}\frac{1}{n}\log\sum_{k=1}^{2^{n}}\left|T_{n,k}\right|^{\theta}
\geq\lim_{n\to\infty}\frac{1}{n}\log\left|T_{n,1}\right|^{\theta}=\lim_{n\to\infty}
\frac{-\theta}{n}\log\left(n+1\right)=0.\]

\textbf{\emph{ad (4)}.} \, In order to determine the left derivative $P^{-}(1)$ of $P$
at $1$,  recall from Proposition \ref{propHMU3} that  $\mu_{\theta}^{*}$ refers to the unique Gibbs measure on
$\Sigma^{*}$
such that
$\mu_{\theta}^{*}\left(C_{n}^{*}\left(y\right)\right)\asymp\exp\left(
-\theta S_{n}I^{*}\left(y\right)-\beta(\theta) S_{n}N^{*}\left(y\right)\right)$,
for all $n \in {\mathbb{N}}, y \in \Sigma^{*}$. For each $n\in\N$, 
let us fix an element $y_X^{(n)}\in\Sigma^{*}$
such that $y^{(n)}_X=\left(X^{n},\ldots\right)$, for $X\in \{A,B\}$. We then 
have by Lemma
\ref{lem:Vergleichbar}, \begin{eqnarray*}
\int N\, d\mu_{\theta}^{*} & = &
\sum_{X\in \{A,B\}}\sum_{n=1}^{\infty} n \cdot \mu_{\theta}^{*}\left(C_{1}^{*}\left(y_X^{(n)}\right)\right)\\
&\asymp&\sum_{n=1}^{\infty}n
\cdot\exp\left(-\theta I^{*}\left(y_A^{(n)}\right)-\beta(\theta) N^{*}\left(y_A^{(n)}\right)\right)\\
   & \gg & \sum_{n=2}^{\infty}n\cdot n^{-2\theta}
   e^{-\beta(\theta)n}\to  \infty, \;\; \mbox{ for }\theta\nearrow 1.\end{eqnarray*}
On the other hand, we have for all $\theta \in (1/2,1]$,
\begin{eqnarray*}
\int I^{*}\, d\mu_{\theta}^{*} & \asymp &
\sum_{X\in \{A,B\}}\sum_{n=1}^{\infty} \log n  \, \mu_{\theta}^{*}\left(C_{1}^{*}\left(y_X^{(n)}\right)\right)\\
&\asymp&\sum_{n=1}^{\infty}\log n\exp\left(-\theta I^{*}\left(y_A^{(n)}\right)-\beta(\theta) N^{*}\left(y_A^{(n)}\right)\right)\\
   & \ll & \sum_{n=1}^{\infty}n^{-2\theta}\log n < \infty.\end{eqnarray*}
   This shows that $P^{-}(1)=0$, and hence $P$ is differentiable
everywhere.

\textbf{\emph{ad (5)}.}  \, Since
$\lim_{\theta\to\infty} P\left(\theta\right)/\theta=0$,
the upper bound of the domain of $\widehat{P}$ is equal to $0$.
For
the lower bound $-\alpha_{+}$ of the domain we have
by \cite{KesseboehmerStratmann:04a} (Proposition 2.3), \begin{equation}
-\alpha_{+}=\lim_{\theta\to-\infty}\frac{P(\theta)}{\theta}=-\sup_{\nu\in\mathcal{M}
\left(\Sigma,\sigma\right)}\int I\, d\nu.\label{eq:alpha+}\end{equation}
 We are left with to determine the actual value of
$\alpha_{+}$. For this, first note that for the linear combination
$m:=1/2
\left(
\delta_{\overline{AB}}+ \delta_{\overline{BA}}\right) \in \mathcal{M}
\left(\Sigma,\sigma\right)$
of the Dirac measures $\delta_{\overline{AB}}$
and $\delta_{\overline{BA}}$ at the periodic points
$\overline{AB}:=\pi^{-1}\left(2-\gamma\right)$ and
$\overline{BA}:=\pi^{-1}\left(\gamma-1\right)$, an elementary
calculation shows that $\int I\, dm=2\log\gamma$.  This implies that
$\sup_{\nu\in\mathcal{M}\left(\Sigma,\sigma\right)}\int I\,
d\nu\geq2\log\gamma$.  For the reverse inequality note that  $\int I\,
d\nu\leq\sup_{x\in\Sigma}\limsup_{n\to\infty}(S_{n}I(x))/n$,
for all
$\nu\in\mathcal{M}\left(\Sigma,\sigma\right)$, where
$\mathcal{M}(\Sigma,\sigma)$ refers to the set of
$\sigma$--invariant Borel probability measures on $\Sigma$.
In order to calculate the right hand side of the latter inequality,
recall that  the smallest interval in 
$\mathcal{T}_{n}$ has the length $(f_{n+1} f_{n+2})^{-1}$.
Using this observation and Corollary
\ref{cor:VergleichTSnI},
we obtain \begin{eqnarray*}
\sup_{y\in\Sigma}\limsup_{n\to\infty}\frac{S_{n}I(y)}{n} &= &\!\!
\sup_{x\in[0,1)}\limsup_{n\to\infty}\frac{-\log\left|T_{n}\left(x\right)\right|}{n}=
\lim_{n\to\infty}
\frac{\log\left(f_{n+1}f_{n+2}\right)}{n}\\
& = &\!\!\!
\lim_{n\to\infty}\frac{\log\left(\gamma^{n+1}-\left(-\gamma\right)^{-(n+1)}\right)+\log
\left(\gamma^{n+2}-\left(-\gamma\right)^{-(n+2)}\right)}{n}\\
  & = & \!\! 2\log\gamma.\end{eqnarray*}
   Note that in here the supremum is achieved at for instance any noble
number in $(0,1)$, that is at numbers whose continued fraction expansion
eventually consists of $1$'s only.

\textbf{\emph{ad (6)}.} \, The result in (3) implies that
   \[
\lim_{\alpha \searrow0}-\widehat{P}(-\alpha)/\alpha=\inf\left\{
t\in\R:P(t)=0\right\} .\]
   Therefore, it is sufficient to show that $1$ is the least zero of $P$. 
For this assume by way of contradiction that $P(s)=0$, for some
$s<1$. Since $P$ is non-increasing, it follows that $P$ vanishes
on the interval $(s,1)$. But this contradicts the fact that $P$ is
real-analytic
on $(-\infty,1)$ and positive at for instance $0$.

\textbf{\emph{ad (7)}.} \,  For all $n\in\N$ 
and $\theta\leq 0$, we have
\begin{eqnarray*}
\left(\frac{\gamma^{n+1}-\left(-\gamma\right)^{-(n+1)}}{\sqrt{5}}\right)^{-2\theta}
& \leq & \left(f_{n+1}f_{n+2}\right)^{-\theta}\leq\sum_{k=1}^{2^{n}}\left|T_{n,k}
\right|^{\theta} \\
& \leq & 2^{n}\left(f_{n+1}f_{n+2}\right)^{-\theta}\leq2^{n}\gamma^{-2\theta
(n+2)}.\end{eqnarray*}
Therefore, \[
-2\theta\log\gamma\leq
P(\theta)\leq\log2-2\theta\log\gamma  \textrm{ for all } \theta\leq0,\]
which implies that $\widehat{P}(-\alpha)\leq0$, for all
$\alpha\in[0,2\log\gamma]$.
Hence, in order to verify that
$\lim_{\alpha\nearrow2\log\gamma}\widehat{P}\left(-\alpha\right)=0$
it is sufficient to show that this limit is non-negative.
For this, let $t(\alpha):=\left(P'\right)^{-1}(-\alpha)$ and recall
that by the variational principle (cf.
\cite{DenkerGrillenbergerSigmund:76}) we have 
for each $\alpha\in[0,2\log\gamma]$ that there exists
$m_{t(\alpha)}\in\mathcal{M}\left(\Sigma,\sigma\right)$
such that \[
P\left(t(\alpha)\right)=h_{m_{t(\alpha)}}-t(\alpha)\int I\,
dm_{t(\alpha)}.\]
Furthermore, by \cite{KesseboehmerStratmann:04a} (Proposition 2.3)
we have $\int I\, dm_{t(\alpha)}=\alpha$. Therefore, if
$\nu\in\mathcal{M}\left(\Sigma,\sigma\right)$
denotes a weak limit of some sequence $\left(\mu_{t(\alpha)}\right)$
for $\alpha$ tending to $2\log\gamma$ from below, then the lower semi-continuity of
the entropy (cf. \cite{DenkerGrillenbergerSigmund:76}) gives \[
h_{\nu}\geq\limsup_{\alpha\nearrow2\log\gamma}h_{m_{t(\alpha)}}=\limsup_{\alpha\nearrow2\log\gamma}
\left(P\left(t(\alpha)\right)+\alpha\cdot
t(\alpha)\right)=\limsup_{\alpha\nearrow2\log\gamma}\left(-\widehat{P}
\left(-\alpha\right)\right).\]
Note that we clearly have $\int I\, d\nu=2\log\gamma$. Now, the final step is to
show that for the discrete measure $m$ considered in the proof of
(5) we have \[
\left\{ \nu\in\mathcal{M}\left(\Sigma,\sigma\right):\int I\,
d\nu=2\log\gamma\right\} =\left\{ m\right\} .\]
This will be sufficient, since $h_{m}=0$.
Therefore, suppose by way of contradiction that there exists $\mu\not=m$
such that \[
\mu\in\left\{ \nu\in\mathcal{M}\left(\Sigma,\sigma\right):\int I\,
   d\nu=2\log\gamma\right\} . \]
 Let us first show that  $\eta:= \mu\left(\left\{ x\in\Sigma:x_{1}=x_{2}=X\right\} \right)>0$,
for some $X \in \{A, B\}$.  If  this would not be the case, 
then the $\sigma$--invariance of 
$\mu$ would imply \[  \mu\left(\left\{ x\in\Sigma:x_{1}=A,x_{2}=B\right\}
\right)= \mu\left(\left\{ x\in\Sigma:x_{1}=B,x_{2}=A\right\} \right)=\frac{1}{2},\] 
and hence we obtain by   induction that $\mu =m$. This contradicts our assumption
$\mu \neq m$, showing that $\eta >0$. We can now continue the 
above argument as follows. 
Since $\{ \nu\in\mathcal{M}\left(\Sigma,\sigma\right):\int I\,
 d\nu=2\log\gamma\} $
is convex, we can assume without loss of generality that $\mu$ is ergodic. 
This then immediately implies that $\lim_{n\to\infty} (S_{n}I(y))/n=\int I\,
d\mu$  for $\mu$-almost every $y\in\Sigma$, and furthermore that for some $X \in \{A, B\}$ and
   $n$ sufficiently large, \begin{equation}
S_{n}\1_{\left\{ x\in\Sigma:x_{1}=x_{2}=X \right\}
}(y)>\frac{n\eta}{2}.\label{eq:ErgodMean}\end{equation}
Let us fix $x\in \Sigma$ with this property, and define  $\tau_k:=\sum _{i=1}^k a_i(x)$ as 
in Lemma \ref{lem:QGammaInequalitystrict}. 
Combining Lemma \ref{lem:TnbyCF} and inequality (\ref{eq:ErgodMean}),
it follows $\left(\tau_k-k-1)\right)\geq \tau_k\eta/2$. Hence, using Corollary
\ref{cor:VergleichTSnI}, Lemma \ref{lem:Vergleichbar} and Lemma
 \ref{lem:QGammaInequalitystrict}, \begin{eqnarray*} 
2\log\gamma & = & \int I\,
d\mu=\lim_{n\to\infty}\frac{S_{n}I(x)}{n}=\lim_{n\to\infty}\frac{-\log
\left|T_{n}(x)\right|}{n}\\
  & = &
\lim_{k\to\infty}\frac{-\log
\left|T_{\tau_k +1}(x)\right|}{\tau_k}= \lim_{k\to\infty}\frac{2\log(q_k(x))}{\tau_k}\\
  & \leq &
 \limsup_{k\to\infty}\frac{2\log\left(\gamma^{\tau_k}\rho^{(\tau_k-k-1)}\right)}{\tau_k}\leq\limsup_{k\to\infty}
\frac{2\log\left(\gamma^{\tau_k}\rho^{\tau_k\eta/2}\right)}{\tau_k}\\
& = & 2\log\gamma+\eta\cdot\log\rho<2\log\gamma.\end{eqnarray*} 
\textbf{\emph{ad (8)}.} \,  First note that $t_{n,2\ell}>t_{n,2\ell\pm1}$,
for each $n\geq2$ and
   $\ell=1,\ldots,2^{n-1}$. This implies that
 $\left|T_{n,2\ell}\right|^{-1}=t_{n,2\ell}\cdot t_{n,2\ell+1}$
and $\left|T_{n,2\ell-1}\right|^{-1}=t_{n,2\ell-1}\cdot t_{n,2\ell}$
are both less than $\left(t_{n,2\ell}\right)^{2}$. Hence, using Lemma
\ref{lem:TnbyCF} and Lemma \ref{lem:QGammaInequalitystrict}, it
follows for $n>2$ and $\theta<0$, \begin{eqnarray*}
\sum_{k=1}^{2^{n}}|T_{n,k}|^{\theta} & \leq &
2\sum_{k=1}^{n}\sum_{A_{k}^{n+1}}q_{k}\left(\left[a_{1},\ldots,a_{k}\right]\right)^{-2\theta}\\
   & \leq & 2\sum_{k=1}^{n}{\binom{n-1}{
k-1}}\left(\gamma^{n+1}\rho^{n+1-k-1}\right)^{-2\theta}\\
   & = & 2\gamma^{-2\theta(n+1)}\sum_{k=0}^{n-1}{\binom{n-1}{
k}}\left(\rho^{n-1-k}\right)^{-2\theta}\\
   & = & 2\gamma^{-2\theta(n+1)}\sum_{k=0}^{n-1}{\binom{n-1}{
k}}\left(\rho^{-2\theta}\right)^{n-1-k}\\
   & \leq &
2\gamma^{-2\theta(n+1)}\left(1+\rho^{-2\theta}\right)^{n-1}.\end{eqnarray*}
Recalling the definition of $P$, we then obtain \[
P\left(\theta\right)\leq-2\theta\log\gamma+\log\left(1+\rho^{-2\theta}\right).\]
  For the lower bound, note that  \[
\sum_{k=1}^{2^{n}}|T_{n,k}|^{\theta}\geq(f_{n+1}f_{n+2})^{-\theta}.\]
   Since $f_{n}=(\gamma^{n}-(-\gamma)^{-n})/\sqrt{5}$, it therefore
follows \[
P\left(\theta\right)\geq-2\theta\log\gamma.\]
   By combining these two bounds for $P\left(\theta\right)$ and then letting $\theta$ tend to
$(-\infty)$,
the proposition follows.
\end{proof}

\section{Multifractal Formalism for continued fractions}

In this section we give the proof of Theorem \ref{Thm:main}, which
we have split up into the three separate parts 
\emph{The lower bound}, \emph{The upper
bound} and \emph{Discussion of boundary points of the spectrum}. We begin with the
following important preliminary remarks.

\begin{rem}\label{rem:Vergleich} 
    $ \, $ \\
    (1) \, 
     Note that by Corollary \ref{cor:VergleichTSnI} and Lemma
\ref{lem:Vergleichbar},
we have for $x\in\Sigma$ and $y\in\Sigma^{*}$ (assuming in each case that the limit 
exists),
\begin{eqnarray*}
\ell_{1}\left(\pi^{*}(y)\right)=\lim_{n\to\infty}\frac{S_{n}I^{*}(y)}{S_{n}N(y)},
&  &
\ell_{2}\left(\pi^{*}(y)\right)=\lim_{n\to\infty}\frac{S_{n}N(y)}{n},\\
\ell_{3}\left(\pi^{*}(y)\right)=\lim_{n\to\infty}\frac{S_{n}I^{*}(y)}{n},
&  &
\ell_{4}\left(\pi(x)\right)=\lim_{n\to\infty}\frac{S_{n}I(x)}{n}.\end{eqnarray*}
(2) \, 
Recall that in \cite{KesseboehmerStratmann:04a} and 
\cite{KesseboehmerStratmann:04} we in particular 
considered oriented geodesics 
$\ell \subset {\mathbb {H}}^{2}$  from $\{\infty\}$ to  [0, 1), and  
coded these by means of their intersections 
with the tesselation given by the $G$-orbit of the fundamental 
domain $F$.  More precisely, if $\ell$ ends at $\xi \in 
[0,1) \cap \I$  such that $\ell$  intersects $ g_{\xi,1}(F), g_{\xi,2}(F), 
g_{\xi,3}(F),\ldots$ in succession, with $g_{\xi,n} \in G$ for all $n 
\in \N$, then $\xi$ is coded by the infinite word $(g_{\xi,1}, 
g_{\xi,2}, 
g_{\xi,3},\ldots)$. Clearly,  this type of coding is 
analogous to the finite coding represented by  $\Sigma$. Hence,
the results of \cite{KesseboehmerStratmann:04a} and 
\cite{KesseboehmerStratmann:04} for the 
Hausdorff dimensions
of the level sets
\[ \mathcal{F}(\alpha):= \left\{ \xi\in [0,1): \frac{d(z_{0}, g_{\xi,n} 
(z_{0}))}{n} = \alpha \right\}\]
can immediately be transfered to the situation in this paper, and in this 
way we obtain that for each
$\alpha\in(0,2\log\gamma)$,
\begin{equation}
\dim_{H}\left(\mathcal{L}_{4}\left(\alpha\right)\right)= \dim_{H}\left( 
\mathcal{F}(\alpha) \right) = 
\frac{\widehat{P}(-\alpha)}{-\alpha}.\label{eq:OldSpec}\end{equation}
Therefore, the following proof of Theorem \ref{Thm:main} will
 in particular also give an alternative proof of the
identity  in (\ref{eq:OldSpec}).
Let us also emphasize that a straight forward inspection of the arguments in the general multifractal analysis 
of \cite{KesseboehmerStratmann:04a}
shows that there we did  not make full use of the group structure of the Kleinian group.
In fact, the arguments  there exclusively
consider certain rooted sub-trees of the Cayley graph of the Kleinian group,
and therefore they continue to hold if the underlying algebraic
structure is only a semi-group acting on hyperbolic space, rather than a group.
Therefore, the main
results of this general multifractal analysis for growth rates can 
 be applied immediately  to the setting in this paper. In this way 
 one also immediately  obtains 
that $P$ is differentiable everywhere, real-analytic on
$(-\infty,1)$
and equal to $0$ otherwise. 
\end{rem}


\subsection{The lower bound}\label{section 5.1}

\begin{lem}
\label{lem:Inclusion} For each $\alpha\in(0,2\log\gamma)$ there
exists a unique Gibbs measure $\mu_{t(\alpha)}^{*}$ on $\Sigma^{*}$
such that for \begin{equation}
\alpha^{*}:=\int I^{*}\,
d\mu_{t(\alpha)}^{*}\,\;\textrm{ and }\;\alpha^{\sharp}:=\int N\,
d\mu_{t(\alpha)}^{*},\label{eq:Alpha*-Alpha}\end{equation}
we have \begin{equation}
\mathcal{L}_{2}\left(\alpha^{\sharp}\right)\cap\mathcal{L}_{3}\left(\alpha^{*}\right)\subset\mathcal{L}_{1}
\left(\alpha\right).\label{eq:SharpSternInclusion}\end{equation}
In here, the function $t$ is given by $t(\alpha):=\left(P'\right)^{-1}(-\alpha)$.
\end{lem}
\begin{proof}
The existence of the unique ergodic Gibbs measure $\mu_{t(\alpha)}^{*}$ has 
already been obtained in Proposition \ref{propHMU2}. 
As shown in Proposition \ref{propHMU3}, the significance of $\mu_{t(\alpha)}^{*}$
is that it allows to represent the Lyapunov exponent $\alpha$ 
in terms of $I^{*}$ and $N$ as follows. 
\begin{equation}
\alpha=-P'(t(\alpha))=\frac{\int I^{*}\, d\mu_{t(\alpha)}^{*}}{\int N\,
d\mu_{t(\alpha)}^{*}}=\frac{\alpha^{*}}{\alpha^{\sharp}}.\label{eq:alpha-alphaalpha}\end{equation}
Using Remark \ref{rem:Vergleich} (1), it follows that
$\mathcal{L}_{2}\left(\alpha^{\sharp}\right)\cap\mathcal{L}_{3}\left(\alpha^{*}\right)\subset\mathcal{L}_{1}\left(\alpha\right)$.

\end{proof}
For the following lemma recall that the Hausdorff dimension
$\dim_{H}\left(\mu\right)$
of a probability measure $\mu$ on some metric space is given by \[
\dim_{H}\left(\mu\right):=\inf\left\{ \dim_{H}(K):\mu(K)=1\right\} .\]
\begin{lem}
For each $\alpha\in(0,2\log\gamma)$ we have, with
$\widetilde{\mu}_{t(\alpha)}:=\mu_{t(\alpha)}^{*}\circ\left(\pi^{*}\right)^{-1}$,
\[
\dim_{H}\left(\widetilde{\mu}_{t(\alpha)}\right)\leq\dim_{H}\left(\mathcal{L}_{2}\left(\alpha^{\sharp}\right)
\cap\mathcal{L}_{3}\left(\alpha^{*}\right)\right)\leq\dim_{H}\left(\mathcal{L}_{1}\left(\alpha\right)\right).\]
\end{lem}
\begin{proof}
The first inequality follows, since by ergodicity of $\mu_{t(\alpha)}^{*}$
we have\[
\widetilde{\mu}_{t(\alpha)}\left(\mathcal{L}_{2}\left(\alpha^{\sharp}\right)\cap
\mathcal{L}_{3}\left(\alpha^{*}\right)\right)=1.\]
The second inequality is an immediate consequence of Lemma
\ref{lem:Inclusion}.
\end{proof}

\begin{lem}
\label{lem:lowerBound} For each $\alpha\in(0,2\log\gamma)$  we have\[
\dim_{H}\left(\widetilde{\mu}_{t(\alpha)}\right)=\frac{\widehat{P}(-\alpha)}{-\alpha}.\]
\end{lem}
\begin{proof}
The aim is to show that the local dimension of 
$\widetilde{\mu}_\alpha$ exists and
is equal to $\widehat{P}(-\alpha)/(-\alpha)$,
for each $\alpha \in(0,2\log\gamma)$.
  For this, let  $B(x,r):=[x-r,x+r]\cap\I$ for $0<r\leq1$ and $x\in\I$,  and
define
\begin{eqnarray*}
m_{r}(x) & := & \max\left\{
n\in\N:\pi^{*}C_{n}^{*}\left(\left(\pi^{*}\right)^{-1}x\right)\supset
B(x,r)\right\} ,\\
n_{r}(x) & := & \min\left\{ n\in\N:\pi^{*}C_{n}^{*}
\left(\left(\pi^{*}\right)^{-1}x\right)\subset B(x,r)\right\}
.\end{eqnarray*}
  We obviously have that $\left|m_{r}(x)-n_{r}(x)\right|$
is uniformly bounded from above, and hence
$\lim_{r\to 0}m_{r}\left(x\right)/n_{r}\left(x\right)=1$.
Combining the Gibbs property of $\mu_{t(\alpha)}^{*}$, Lemma
\ref{lem:Vergleichbar},  (\ref{eq:Alpha*-Alpha})
and (\ref{eq:alpha-alphaalpha}), it follows
   for $\widetilde{\mu}_{t(\alpha)}$-almost every $x$,
\\
\\
${\displaystyle
\limsup_{r\to 0}\frac{\log\widetilde{\mu}_{t(\alpha)}\left(B(x,r)\right)}{\log
r}}$\begin{eqnarray*}
   & \leq &
   \limsup_{r\to 0}\frac{-t(\alpha)\left(S_{n_{r}\left(x\right)} I^{*}\left(\left(\pi^{*}\right)^{-1}x\right)
   \right)-P\left(t\left(\alpha\right)\right)S_{n_{r}\left(x\right)} N\left(\left(\pi^{*}\right)^{-1}x\right)}
   {-\left(S_{m_{r}\left(x\right)}I^{*}(x)\right)}\\
  & = &
\limsup_{r\to 0}\frac{-t(\alpha)\frac{S_{n_{r}\left(x\right)}I^{*}\left(\left(\pi^{*}\right)^{-1}x\right)}
{S_{n_{r}\left(x\right)}N\left(\left(\pi^{*}\right)^{-1}x\right)}-P\left(t\left(\alpha\right)\right)}
{-\frac{S_{n_{r\left(x\right)}}I^{*}\left(\left(\pi^{*}\right)^{-1}x\right)}{S_{n_{r}\left(x\right)}N
\left(\left(\pi^{*}\right)^{-1}x\right)}\cdot\frac{S_{m_{r}\left(x\right)}I^{*}\left(\left(\pi^{*}
\right)^{-1}x\right)}{m_{r}(x)}\frac{n_{r}(x)}{S_{n_{r}\left(x\right)}
N\left(\left(\pi^{*}\right)^{-1}x\right)}\cdot\frac{m_{r}\left(x\right)}{n_{r}\left(x\right)}}\\
  & = &
\frac{t(\alpha)\alpha+P(t(\alpha))}{\alpha}=\frac{\widehat{P}(-\alpha)}{-\alpha}
.\end{eqnarray*}
The reverse inequality for the `$\liminf$' is obtained by 
similar means, and we omit its proof.  
\end{proof}

\subsection{The upper bound}\label{section 5.2}

\begin{lem}
\label{lem:UpperBound} For each $\alpha\in(0,2\log\gamma)$ we
have\begin{eqnarray*}
\dim_{H}\left(\pi^{*}\left\{
x\in\Sigma^{*}:\liminf_{n\to\infty}\frac{S_{n}I^{*}(x)}{S_{n}N(x)}\geq\alpha\right\}
\right) & \leq & \frac{\widehat{P}(-\alpha)}{-\alpha}.\end{eqnarray*}
\end{lem}
\begin{proof}
Note that $\max\left\{
t\left(\alpha\right)+P(t\left(\alpha\right))/s:s\in[\alpha,2\log\gamma)\right\}
=t(\alpha)+P(t(\alpha))/\alpha$,
for each $\alpha\in(0,2\log\gamma)$. By combining this with 
the Gibbs property of $\mu_{t(\alpha)}^{*}$,
it follows that for each $\varepsilon>0$ and $x\in\Sigma^{*}$ such
that
$\pi^{*}(x)\in\mathcal{L}_{4}\left(\alpha\right)$, we have for $n$ 
sufficiently large, \begin{eqnarray*}
\mu_{t(\alpha)}^{*}\left(C_{n}^{*}(x)\right) & \gg &
\exp\left(-t\left(\alpha\right)\, S_{n}I^{*}(x)-P(t(\alpha))\,
S_{n}N(x)\right)\\
   & = &
\exp\left(-S_{n}I^{*}(x)\left(t\left(\alpha\right)+P(t\left(\alpha\right))\frac{S_{n}N(x)}{S_{n}I^{*}(x)}\right)\right)\\
  & \gg &
\left(\exp\left(-S_{n}I^{*}(x)\right)\right)^{\frac{\widehat{P}(-\alpha)}{-\alpha}+\varepsilon}\\
  & \gg &
\left|\pi^{*}\left(C_{n}^{*}(x)\right)\right|^{\frac{\widehat{P}(-\alpha)}{-\alpha}+\varepsilon}.\end{eqnarray*}
Therefore,  for the sequence of balls 
$\left(B\left(\pi(x),r_{n}\right)\right)$ 
with radii $r_{n}:=\left|\pi^{*}\left(C_{n}^{*}(x)\right)\right|$ and 
centre $\pi(x)$, which tends to $\{\pi(x)\}$ for $n$ tending to infinity,
we have \[ \widetilde{\mu}_{t(\alpha)}\left(
B\left(\pi(x),r_{n}\right)\right)\gg\mu_{t(\alpha)}^{*}
\left(C_{n}^{*}(x)\right)\gg\left(r_{n}\right)^{\frac{\widehat{P}
(-\alpha)}{-\alpha}+\varepsilon}.\] Applying the mass distribution
principle, the proposition follows.
\end{proof}
\begin{cor}
\label{cor:DimUpperBound} For each $\alpha\in(0,2\log\gamma)$ we
have \[
\max\left\{
\dim_{H}\left(\mathcal{L}_{2}\left(\alpha^{\sharp}\right)\cap\mathcal{L}_{3}\left(\alpha^{*}\right)\right),\dim_{H}\left(\mathcal{L}_{1}\left(\alpha\right)\right)\right\}
\leq\frac{\widehat{P}(-\alpha)}{-\alpha}.\]
\end{cor}
\begin{proof}
The assertion is an immediate consequence of combining Lemma \ref{lem:UpperBound}
and the fact \[
\mathcal{L}_{2}\left(\alpha^{\sharp}\right)\cap\mathcal{L}_{3}\left(\alpha^{*}\right)\subset\mathcal{L}_{1}(\alpha)\subset\left\{
x\in\Sigma^{*}:\liminf_{n\to\infty}\frac{S_{n}I^{*}(x)}{S_{n}N(x)}\geq\alpha\right\}.\]
\end{proof}

\subsection{Discussion of  boundary points of the
spectrum.\label{sub:Discussion-of-boundary-points}}
$ \, $\\
 {\em The case $\alpha=0$}: \,  Recall the two classical
results of L{\'e}vy and Khintchin mentioned in the introduction. 
>From these we immediately deduce that $\tau(0)=1$. Also,
recall that by Proposition \ref{pro:AnalyticPropertiesP}
(6) we have that
$\lim_{\alpha\searrow0}\widehat{P}(-\alpha)/(-\alpha)=1$.
This shows that
$\tau(0)=\lim_{\alpha\searrow0} \widehat{P}(-\alpha)/(-\alpha) =1$,
and hence gives  that the dimension function $\tau$ is continuous from
the right at $0$.\\
In order to show that $\alpha^{*}\left(0\right)=\chi$, we argue as 
follows. For $\alpha=0$, we already know that $0=\int I^{*}\, d\mu_{1}^{*}/\int N\,
d\mu_{1}^{*}= \alpha^{*}\left(0\right)/\infty$ and that
$\lim_{k\to\infty} (2\log q_{k}(x)) / k=\alpha^{*}(0)$, for
$\mu_{1}^{*}\circ\left(\pi^{*}\right)^{-1}$-almost every
$x\in\left(0,1\right)$.  Hence,  L{\'e}vy's result gives that,
if $\mu_{1}^{*}\circ\left(\pi^{*}\right)^{-1}$ is absolutely 
continuous with respect
to the Lebesgue measure $\lambda$ on $(0,1)$ then  
$\alpha^{*}\left(0\right)=\chi$.  Hence, it remains to show that $\mu_{1}^{*}\circ\left(\pi^{*}\right)^{-1}$
has this property. For this, consider some $T\in\mathcal{T}_{n}$ for 
$n\in\N$, and fix $y\in\Sigma^{*}$ and $k\in\N$
such that $\pi^{*}\left(C_{k}^{*}(y)\right)=T\cap\I$. Using the Gibbs property of $\mu_{1}^{*}$
and
Lemma \ref{lem:Vergleichbar}, we  obtain \begin{eqnarray*}
\mu_{1}^{*}\circ\left(\pi^{*}\right)^{-1}\left(T\right) & \asymp &
\mu_{1}^{*}\left(C_{k}^{*}\left(y\right)\right)\asymp\exp\left(-S_{k}\left(I^{*}(y)\right)\right)\\
& \asymp &
\left|\pi^{*}\left(C_{k}^{*}\left(y\right)\right)\right|\asymp\lambda\left(T\right).
\end{eqnarray*}
{\em The case $\alpha=2\log\gamma$}: \,  In order to show that the dimension
function $\tau$ is continuous from the left at $2\log\gamma$, we proceed as follows. Proposition
\ref{pro:AnalyticPropertiesP}
(7) implies that
$\lim_{\alpha\nearrow2\log\gamma}\widehat{P}(-\alpha)/(-\alpha)=0$.
  Using monotonicity of the Hausdorff 
  dimension together with Lemma \ref{lem:UpperBound},
 it then follows  \[
0\leq\tau\left(2\log\gamma\right)\leq\lim_{\alpha\nearrow2\log\gamma}\tau(\alpha)=0.\]
Hence,  we have $\tau(2\log\gamma)=0$, which gives that $\tau$ is continuous 
from the left at $2\log\gamma$.\\
Finally, for the left derivative of $\tau$ at $2\log\gamma$, 
note that a straight forward computation of
 the derivative of $\tau$ on the interval $(0,2\log\gamma)$ shows that $\tau'(\alpha)=-P\left(t\left(\alpha\right)\right)/\alpha^{2}$.
Since $t(\alpha)$ tends to $(-\infty)$ as $\alpha$ approaches $2\log\gamma$, it
follows 
\[\lim_{\alpha\nearrow2\log\gamma}\tau'(\alpha)=-\infty.\]

\section{Multifractal formalism for approximants}

In this section we outline the necessary changes which have to be 
implemented in the proof  of Theorem \ref{Thm:main} in order to 
derive Theorem \ref{thm:main3}.

The analytic properties of $P_{D}$ as stated in
Theorem \ref{thm:main3} can be obtained as follows. In 
Section 4 replace the function 
$N:\overline{\Sigma}\to\N$ (resp. $N:\Sigma^{*} \to \N$) by the
function $\overline{{\bf 1}}: \overline{\Sigma}\to \{1\}$ (resp. 
${\bf 1}^{*}: \Sigma^{*} \to \{1\}$)  constant equal to $1$. 
In this way we 
obtain for the pressure 
function $\overline{\P}$ 
associated
with $\overline{\Sigma}$,

\[
\overline{\P}\left(-\theta \overline{I}-\overline{\P}
\left(-\theta
\overline{I}\right) \overline{{\bf 1}}\right)=0.\]
Also, note that by Lemma \ref{lem:Vergleichbar} we have \[
P_{D}(\theta)=\overline{\P} \left(-\theta \overline{I}\right).\]
(Below, we shall specify the domain  of $P_D$).  Hence, combining these observations with Proposition 
\ref{propHMU3} adapted to the situation here,
the alleged  analytic properties of $P_{D}$
follow. Also, using the same strategy in Section \ref{section 5.1} and 
\ref{section 5.2}, that is replacing in there the function $N$ by 
the function ${\bf 
1}^{*}$, one immediately obtains 
  \[
\tau_{D}(\alpha)=
\frac{\widehat{P}_{D}(-\alpha)}{-\alpha}.\] 
(Below, we shall specify the domain  of $\tau_D$).

For clarifying the range of $P_D$ and of $\tau_D$, and for the discussion of the boundary points of $\tau_D$,
we first remark that $P_{D}$ has a singularity at $1/2$.
This follows, since for every approximant
$\left[a_{1},\ldots,a_{k}\right]$
we have (see e.g. \cite{Khihtchine:36}) \[ \prod_{i=1}^{k}a_{i}\leq
q_{k}\left(\left[a_{1},\ldots,a_{k}\right]\right)\leq2^{k}\prod_{i=1}^{k}a_{i},\]
which immediately gives 
\[
0\leq\log\zeta\left(\theta\right)-P_{D}\left(\theta\right)\leq2\theta\log2,
\textrm{ for }\theta>1/2.\] Here, $\zeta$ refers to the Riemann 
zeta-function $\zeta(\theta):= \sum_{n \in \N} n^{-2\theta}$. This shows that $P_{D}(\theta)$ and 
$P_{D}'(\theta)$  both tend to infinity for $\theta$ tending to $1/2$ from 
above.  From this we deduce that 
$\widehat{P}_{D}\left(-\alpha\right)$ is well defined for arbitrary
large values of $\alpha$, and also that
\[\lim_{\alpha \to \infty}\widehat{P}_{D}(-\alpha)/(-\alpha)=1/2.\] 

In order to see that the domain of $\widehat{P}_{D}$ is the interval
$[2\log\gamma,\infty)$ and that $\lim_{\alpha\searrow
2\log\gamma}\widehat{P}_{D}(-\alpha)/(-\alpha)=0$, it is now sufficient to
show that \begin{equation}
\lim_{\theta\to\infty}\left|P_{D}\left(\theta\right)+2
\theta\log\gamma\right|=0.
\label{eq:toshow}\end{equation}
Indeed, on the one hand we have \begin{eqnarray*}
\lim_{k\to\infty}\frac{1}{k}\log\sum_{\left[a_{1},\ldots,a_{k}\right]}q_{k}
\left(\left[a_{1},\ldots,a_{k}\right]\right)^{-2\theta}
& \leq & \lim_{k\to\infty}-\frac{1}{k}2\theta\log
q_{k}\left(\gamma\right)=-2\theta\log\gamma.\end{eqnarray*} On the
other hand, using Lemma \ref{lem:QGammaInequalitystrict} and 
\ref{lem:TnbyCF}, we observe for $N\in \N$ and  $\theta>(1+\log
N)/(2\log\gamma)$, \\
$\quad{\displaystyle
\lim_{k\to\infty}\frac{1}{k}\log\sum_{\left[a_{1}
,\ldots,a_{k}\right]}q_{k}\left(\left[a_{1},\ldots
,a_{k}\right]\right)^{-2\theta}}$\begin{eqnarray*}
\qquad\qquad & \leq &
\limsup_{k\to\infty}\frac{1}{k}\log\sum_{n=k+1}^{\infty}\binom{n}{k}\gamma^{-2\theta n}\\
   & = &
-2\theta\log\gamma+\limsup_{k\to\infty}\frac{1}{k}\log\sum_{n=1}^{\infty}\binom{n+k}{k}
\gamma^{-2\theta n}\\
   & \leq &
-2\theta\log\gamma+\limsup_{k\to\infty}\frac{1}{k}\log\sum_{n=1}^{\infty}\frac{\left(n+k\right)^{\left(n+k\right)}}{k^{k}n^{n}}\gamma^{-2\theta
n}\\
   & \leq &
-2\theta\log\gamma+\limsup_{k\to\infty}\frac{1}{k}\log\sum_{n=1}^{\infty}\left(1+\frac{k}{Nn}\right)^{n}N^{n}\left(1+\frac{n}{k}\right)^{k}\gamma^{-2\theta
n}\\
   & \leq 
&-2\theta\log\gamma+\limsup_{k\to\infty}\frac{1}{k}\log\sum_{n=0}^{\infty}e^{k/N}e^{n\left(1+\log
N-2\theta\log\gamma\right)}\\
   & \leq & -2\theta\log\gamma+1/N.\end{eqnarray*}
A combination of these two observations gives the statement in (\ref{eq:toshow}).

In order to prove continuity of $\tau_{D}$ at $2\log\gamma$,
note that by arguing similar as in the proof of  Lemma
\ref{lem:UpperBound}, we obtain
 for $\alpha<\chi$, \begin{eqnarray*}
\dim_{H}\left(\pi_{\mathrm{CF}} \left\{
x\in\overline{\Sigma}:\limsup_{n\to\infty}\frac{S_{n}\overline{I}(x)}{n}\leq\alpha\right\}
\right) & \leq &
\frac{\widehat{P}_{D}(-\alpha)}{-\alpha}.\end{eqnarray*}
Combining this  with the monotonicity of Hausdorff
dimension, it follows 
\[ \dim_{H}\left(\mathcal{L}_{3}
\left(2\log\gamma\right)\right)\leq\lim_{\alpha\searrow2\log\gamma}
  \widehat{P}_{D}(-\alpha)/(-\alpha)=0.\]

Finally,   the same argument as used  in
Section \ref{sub:Discussion-of-boundary-points} for  
determining the limit behaviour of $\tau'$, gives that for the 
left derivative of $\tau_{D}$ at $2\log\gamma$ we have  
\[ \lim_{\alpha\searrow2\log\gamma}\tau_{D}' 
\left(\alpha\right)=\infty.\]
This finishes the proof of Theorem \ref{thm:main3}.


\newcommand{\etalchar}[1]{$^{#1}$}

\end{document}